\documentclass[11pt]{article}

\usepackage{amsmath,amsfonts,amsthm,amssymb,graphicx,tikz}
\usepackage[all,2cell,ps]{xy}

\theoremstyle{plain}
\newtheorem{thm}{Theorem}[section]

\newtheorem{lem}[thm]{Lemma}
\newtheorem{prop}[thm]{Proposition}

\newtheorem{cor}[thm]{Corollary}
\newtheorem*{thma}{Theorem A}
\newtheorem*{thmb}{Theorem B}
\newtheorem*{thmc}{Theorem C}

\newtheorem{proj}{Project}
\newtheorem{qtn}[proj]{Question}

\theoremstyle{definition}

\newtheorem{rem}[thm]{Remark}

\theoremstyle{remark}

\newcommand{\bbF}{\mathbb{F}}

\newcommand{\bbP}{\mathbb{P}}
\newcommand{\bbQ}{\mathbb{Q}}
\newcommand{\bbR}{\mathbb{R}}

\newcommand{\bbZ}{\mathbb{Z}}

\newcommand{\calB}{\mathcal{B}}
\newcommand{\calC}{\mathcal{C}}

\newcommand{\calJ}{\mathcal{J}}

\newcommand{\calL}{\mathcal{L}}

\newcommand{\calO}{\mathcal{O}}
\newcommand{\calP}{\mathcal{P}}
\newcommand{\calQ}{\mathcal{Q}}
\newcommand{\calR}{\mathcal{R}}

\newcommand{\frakA}{\mathfrak{A}}

\newcommand{\frakS}{\mathfrak{S}}

\newcommand{\frakh}{\mathfrak{h}}

\newcommand{\frakp}{\mathfrak{p}}

\newcommand{\fraksl}{\mathfrak{sl}}

\newcommand{\al}{\alpha}
\newcommand{\gam}{\gamma}
\newcommand{\Gam}{\Gamma}
\newcommand{\del}{\delta}

\newcommand{\Del}{\Delta}
\newcommand{\ep}{\epsilon}

\newcommand{\Lam}{\Lambda}
\newcommand{\sig}{\sigma}

\DeclareMathOperator{\M}{M}
\DeclareMathOperator{\SL}{SL}
\DeclareMathOperator{\PSL}{PSL}

\DeclareMathOperator{\Id}{Id}

\DeclareMathOperator{\Tr}{Tr}
\DeclareMathOperator{\tr}{tr}

\newcommand{\conj}{\overline}
\newcommand{\wh}{\widehat}
\newcommand{\wt}{\widetilde}

\newenvironment{pf}{\begin{proof}}{\end{proof}}

\newenvironment{enum}{\begin{enumerate}}{\end{enumerate}}

\usepackage{tikz-cd}
\allowdisplaybreaks

\makeatletter
\let\@@pmod\pmod
\DeclareRobustCommand{\pmod}{\@ifstar\@pmods\@@pmod}
\def\@pmods#1{\mkern4mu({\operator@font mod}\mkern 6mu#1)}
\makeatother

\usetikzlibrary{arrows}

\title{Geometry of the Wiman--Edge monodromy}
\author{Matthew Stover\footnote{This material is based upon work supported by Grant Number DMS-1906088 from the National Science Foundation.} \\ \small{Temple University}\\ \small{\textsf{mstover@temple.edu}}}
\date{\today}

\begin{document}

\maketitle

%%%%%%%%%%%%%%%%%%%%
\begin{abstract}
The Wiman--Edge pencil is a pencil of genus $6$ curves for which the generic member has automorphism group the alternating group $\frakA_5$. There is a unique smooth member, the \emph{Wiman sextic}, with automorphism group the symmetric group $\frakS_5$. Farb and Looijenga proved that the monodromy of the Wiman--Edge pencil is commensurable with the Hilbert modular group $\SL_2(\bbZ[\sqrt{5}])$. In this note, we give a complete description of the monodromy by congruence conditions modulo $4$ and $5$. The congruence condition modulo $4$ is new, and this answers a question of Farb--Looijenga. We also show that the smooth resolution of the Baily--Borel compactification of the locally symmetric manifold associated with the monodromy is a projective surface of general type. Lastly, we give new information about the image of the period map for the pencil.
\end{abstract}
%%%%%%%%%%%%%%%%%%%%

%%%%%%%%%%%%%%%%%%%%
\section{Introduction}\label{sec:Intro}
%%%%%%%%%%%%%%%%%%%%

The Wiman--Edge pencil is a pencil of genus $6$ curves for which the generic member has automorphism group the alternating group $\frakA_5$. There is a unique smooth member, the \emph{Wiman sextic}, whose automorphism group is the symmetric group $\frakS_5$. In a recent series of papers, Dolgachev, Farb, and Looijenga \cite{DFL} then Farb and Looijenga \cite{FL1, FL2} studied this pencil from a modern perspective. See these papers for basic facts about the pencil.

In \cite{FL1}, Farb and Looijenga proved that the monodromy of the Wiman--Edge pencil is commensurable with the Hilbert modular group $\SL_2(\bbZ[\sqrt{5}])$. In this paper, we give a complete description of the monodromy by congruence conditions modulo $4$ and $5$. The congruence condition modulo $4$ is new, and this answers a question of Farb--Looijenga \cite[Qn.~4.4]{FL1}.

Before stating our main result, we need some notation. Let $\calO$ be the ring of integers in $\bbQ(\sqrt{5})$ and $\calO_o$ denote the index two subring $\bbZ[\sqrt{5}]$. There is a unique prime ideal $\frakp_5$ of $\calO$ with residue field $\bbF_5$, the field with five elements. The following is a precise description of the Wiman--Edge monodromy.

%%%%%%%%%%%%%%%%%%%%
\begin{thma}\label{thm:DescriptionIntro}
The monodromy group $\Gam$ of the Wiman--Edge pencil can described equivalently in $\SL_2(\calO)$ and $\SL_2(\calO_o)$ by the following two properties:
\begin{enumerate}

\item The reduction of $\Gam$ modulo $\frakp_5$ is unipotent.

\item The reduction $\conj{\Gam}$ of $\Gam$ modulo $4 \calO$ fits into an exact sequence
\[
1 \to C_4 \to \conj{\Gam} \to \SL_2(\bbF_2) \to 1
\]
where $C_4 < \fraksl_2(\bbF_4)$ denotes the subgroup of matrices $\left(\begin{smallmatrix} x_1 & x_2 \\ x_3 & x_1 \end{smallmatrix}\right) \in \fraksl_2(\bbF_4)$ so that $\Tr_{\bbF_4 / \bbF_2}(x_1 + x_2 + x_3) = 0$, where $\Tr_{\bbF_4 / \bbF_2}$ denotes the trace.

\end{enumerate}
\end{thma}
%%%%%%%%%%%%%%%%%%%%a

The proof of Theorem~A is in \S\ref{sec:Mono}. Other than the explicit matrix generators for $\Gam$ given in \cite{FL1}, our proof that the monodromy has finite index in $\SL_2(\calO)$ is completely independent of Farb and Looijenga's argument. We use a presentation for $\SL_2(\calO)$ due to Yoshida \cite{Yoshida} along with the computer algebra program Magma \cite{Magma} to describe $\Gam$. It would be interesting to understand the relationship between the condition on $\Gam$ modulo $4$ and the geometry of the pencil.

%%%%%%%%%%%%%%%%%%%%
\begin{qtn}
Can one interpret the mod $4$ condition on $\Gam$ in terms of properties of the Wiman--Edge pencil analogous to the description of the mod $\frakp_5$ condition described in \cite[Cor.~4.3]{FL1}?
\end{qtn}
%%%%%%%%%%%%%%%%%%%%

Our methods allow us to study other properties of the monodromy. In \S\ref{sec:Quo}-\ref{sec:Resolve}, we study the smooth resolution $Y_\Gam$ of the Baily--Borel compactification of the locally symmetric manifold $X_\Gam$ associated with the monodromy. We prove:

%%%%%%%%%%%%%%%%%%%%
\begin{thmb}\label{thm:ChernNosIntro}
The smooth compactification $Y_\Gam$ of $X_\Gam$ is a smooth projective surface of general type with Chern numbers:
\begin{align*}
c_1^2(Y_\Gam) &= 16 \\
c_2(\Gam) &= 56
\end{align*}
It has holomorphic Euler characteristic $\chi(\calO_{Y_\Gam}) = 6$, irregularity $q(Y_\Gam) = 0$, and geometric genus $p_g(Y_\Gam) = 5$.
\end{thmb}
%%%%%%%%%%%%%%%%%%%%

Finally, in \S\ref{sec:Image}, we study the period map of the monodromy, following the suggestion of Farb--Looijenga that one should understand the image of the period map in terms of the \emph{Klein plane} \cite[\S 5]{FL1}. Briefly, the period map can be considered as a closed embedding of $\bbP^1$ into the Baily--Borel compactification $Y^\circ$ of $(\frakh \times \frakh) / \SL_2(\calO_o)$. If $X[2]$ denotes the quotient of $\frakh \times \frakh$ by the level two congruence subgroup of $\SL_2(\calO)$ and $Y[2]$ its smooth compactification, then there is an action of $S_3 = \SL_2(\bbF_2)$ on $Y[2]$ inducing a rational map to $Y^\circ$. We will show:

%%%%%%%%%%%%%%%%%%%%
\begin{thmc}\label{thm:IntroPeriod}
The image of the period map of the Wiman--Edge pencil is determined by a unique $S_3$-invariant curve of genus one on $Y[2]$.
\end{thmc}
%%%%%%%%%%%%%%%%%%%%

See Proposition~\ref{prop:Bcover} for the precise statement. Hirzebruch famously identified $Y[2]$ with the blowup of the Clebsch cubic surface at its ten Eckardt points (e.g., see \cite[\S VIII.2]{vdG}). The Klein plane is $\bbP^2$ realizing $Y[2]$ as the blowup of the plane at the sixteen points associated with an icosahedron (where the blowup at these points is the Clebsch cubic) and a dually inscribed dodecahedron (i.e., the points that transform to the Eckardt points on the Clebsch cubic). In particular, we see that the image of the period map is associated with a unique $S_3$-invariant plane curve of genus one, which we call the \emph{Wiman--Edge monodromy plane curve}.

We end the paper by proving a number of other general properties of the monodromy plane curve, in particular regarding the combinatorics of how it intersects the classical points and lines on $\bbP^2$ associated with the icosahedron and its dual dodecahedron. We were not able to completely describe the curve, most critically because we do not know if the monodromy plane curve is smooth. If it is smooth, Proposition~\ref{prop:Image1} proves that it must be one of two smooth plane cubic curves passing through the unique $S_3$-orbit of vertices of the dodecahedron with cardinality six. It would of course be interesting to find further conditions that would precisely describe the monodromy plane curve, and hence completely describe the image of the period map.

\medskip

We close the introduction with a question related to a curious discovery made while writing this paper. The Wiman curve is naturally uniformized as a quotient of the hyperbolic plane by a subgroup of the $(2,4,6)$ triangle group (e.g., see \cite[\S 4.3]{FL2}). This implies that the Wiman curve is uniformized by a cocompact arithmetic Fuchsian group with associated quaternion algebra the unique quaternion algebra $A$ over $\bbQ$ ramified at the primes $\{2, 3\}$ (see \cite[\S 13.3]{MaclachlanReid}). Classical embedding results for quaternion algebras imply that $A$ embeds as a subalgebra of $\M_2(k)$, since $2$ and $3$ are inert in $k$.

Specifically, there is a maximal order $\calO$ of $A$ and a subgroup $\Lam$ of the group $N(\calO^1)$, the normalizer in $\SL_2(\bbR)$ of the group $\calO^1$ of units of norm one in $\calO$, so that $C_0 = \frakh / \Lam$ is the Wiman curve. Let $\Lam^1$ be the index two subgroup of $\Lam$ given by intersecting $\Lam$ with $\calO^1$. We found an embedding $A \hookrightarrow \M_2(k)$ so that $\Lam^1$ maps into $\SL_2(\calO)$. However, the image is not contained in $\SL_2(\calO_o)$, so this does not realize $\Lam$ as a subgroup of the monodromy group $\Gam$ of the Wiman--Edge pencil, only some finite index subgroup of $\Lam$. It follows that there is a finite \'etale cover $C^\prime \to C_0$ and a closed immersion
\[
C^\prime \looparrowright X_\Gam = (\frakh \times \frakh) / \Gam
\]
with image a totally geodesic submanifold of $X_\Gam$. It would be interesting to know if this curve $C^\prime$ has any significance for the pencil or the image of the period map.

%%%%%%%%%%%%%%%%%%%%
\begin{qtn}
Let $D^\prime \subset X_\Gam$ be an immersed totally geodesic curve commensurable with the Wiman curve. Does the arithmetic Fuchsian subgroup of the monodromy $\Gam$ of the Wiman--Edge pencil associated with $D^\prime$ have special significance? What about the points of intersection of $D^\prime$ with the image of the period map for the pencil?
\end{qtn}
%%%%%%%%%%%%%%%%%%%%

%%%%%%%%%%%%%%%%%%%%
\subsubsection*{Acknowledgments}
I thank Benson Farb and Eduard Looijenga for a number of interesting conversations about the Wiman--Edge pencil, and for their encouragement to write up the results in this paper.
%%%%%%%%%%%%%%%%%%%%

%%%%%%%%%%%%%%%%%%%%
\section{Basic facts and notation}\label{sec:Facts}
%%%%%%%%%%%%%%%%%%%%

We follow the basic notation in \cite{FL1}, which we now recall and fix for the remainder of the paper. Set $k = \bbQ(X)$, where $X^2 - X - 1 = 0$. Then the ring of integers $\calO$ of $k$ is $\bbZ[X]$, and the unit group $\calO^*$ is isomorphic to $\bbZ \times \bbZ / 2$ generated by $X$ and $-1$, where $X^{-1} = X - 1$.

In what follows $\bbF_q$ will always denote the field with $q$ elements. Two prime ideals of $\calO$ that will appear throughout this paper are $2 \calO$, which has residue field $\bbF_4$, and $\frakp_5 = (1 - 2 X) \calO$, which has residue field $\bbF_5$. Note that $(1 - 2X)^2 = 5$ which explicitly realizes the isomorphism $k \cong \bbQ(\sqrt{5})$.

Let $\calO_o$ denote the additive subgroup $\bbZ[1, 2 X] = \bbZ[1, \sqrt{5}]$ and observe that $2\calO \subset \calO_o$. The isomorphism $\bbF_4 \cong \bbF_2[X]$ characterizes $\calO_o$ as the pullback of $\bbF_2 \subset \bbF_4$ under reduction modulo $2 \calO$. This interpretation will be used often in what follows. Notice moreover that $\calO_o$ maps onto $\calO / \frakp_5$, since it maps to a nontrivial cyclic additive subgroup of $\bbF_5$.

\medskip

We now consider $\SL_2(\calO_o) < \SL_2(\calO)$, and begin by recording the following basic lemma that follows from the above characterization of $\calO_o$ as a subgroup of $\calO$.

%%%%%%%%%%%%%%%%%%%%
\begin{lem}\label{lem:FindOo}
The reduction homomorphism
\[
h_2 : \SL_2(\calO) \to \SL_2(\calO / 2 \calO) \cong \SL_2(\bbF_4)
\]
is surjective and $\SL_2(\calO_o)$ is the pullback of the subgroup $\SL_2(\bbF_2) < \SL_2(\bbF_4)$ induced by the field inclusion $\bbF_2 \hookrightarrow \bbF_4$.
\end{lem}
%%%%%%%%%%%%%%%%%%%%

Following Yoshida \cite{Yoshida}, we now give a presentation for $\PSL_2(\calO)$ and use it to derive  one for $\SL_2(\calO)$. First, consider the following matrices in $\SL_2(\calO)$:
\begin{align*}
\sig &= \begin{pmatrix} 0 & 1 \\ -1 & 0 \end{pmatrix} & \tau &= \begin{pmatrix} 1 & 1 \\ 0 & 1 \end{pmatrix} \\
\mu &= \begin{pmatrix} X & 0 \\ 0 & X - 1 \end{pmatrix} & \eta &= \begin{pmatrix} 1 & X \\ 0 & 1 \end{pmatrix}
\end{align*}
Notice that $\langle \sig, \tau \rangle = \SL_2(\bbZ)$. Let $\wh{g}$ denote the image in $\PSL_2$ of a given $g \in \SL_2$ and $[\, ,\, ]$ denote the commutator. Yoshida showed:

%%%%%%%%%%%%%%%%%%%%
\begin{thm}[Thm.~5.1 \cite{Yoshida}]\label{thm:PresentPSL}
The group $\PSL_2(\calO)$ is generated by $\wh{\sig}$, $\wh{\mu}$, $\wh{\tau}$, and $\wh{\eta}$ subject to the relations:
\begin{align*}
\wh{R}_1 &= \wh{\sig}^2 & \wh{R}_2 &= (\wh{\sig} \wh{\tau})^3 \\
\wh{R}_3 &= (\wh{\sig} \wh{\mu})^2 & \wh{R}_4 &= [\wh{\tau}, \wh{\eta}] \\
\wh{R}_5 &= \wh{\mu} \wh{\tau} \wh{\mu}^{-1} (\wh{\tau} \wh{\eta})^{-1} & \wh{R}_6 &= \wh{\mu} \wh{\eta} \wh{\mu}^{-1} (\wh{\tau} \wh{\eta}^2)^{-1} \\
\wh{R}_7 &= \wh{\sig} \wh{\eta} \wh{\sig} (\wh{\tau} \wh{\eta}^{-1} \wh{\sig} \wh{\eta}^{-1} \wh{\mu})^{-1} &  &
\end{align*}
\end{thm}
%%%%%%%%%%%%%%%%%%%%

We now use Theorem~\ref{thm:PresentPSL} to give a presentation for $\SL_2(\calO)$ using the central element
\[
z_0 = \begin{pmatrix} -1 & 0 \\ 0 & -1 \end{pmatrix}
\]
of order two. We then have the following.

%%%%%%%%%%%%%%%%%%%%
\begin{cor}\label{cor:PresentSL}
The group $\SL_2(\calO)$ is generated by $z_0$, $\sig$, $\mu$, $\tau$, and $\eta$ subject to the relations:
\begin{align*}
C_0 &= z_0^2 & & \\
C_1 &= [z_0, \sig] & C_2 &= [z_0, \mu] \\
C_3 &= [z_0, \tau] & C_4 &= [z_0, \eta] \\
R_1 &= \sig^2 z_0 & R_2 &= (\sig \tau)^3 \\
R_3 &= (\sig \mu)^2 z_0 & R_4 &= [\tau, \eta] \\
R_5 &= \mu \tau \mu^{-1} (\tau \eta)^{-1} & R_6 &= \mu \eta \mu^{-1} (\tau \eta^2)^{-1} \\
R_7 &= \sig \eta \sig (\tau \eta^{-1} \sig \eta^{-1} \mu)^{-1} z_0 &  &
\end{align*}
\end{cor}
%%%%%%%%%%%%%%%%%%%%

%%%%%%%%%%%%%%%%%%%%
\begin{pf}
The relations indicating that $z_0$ defines a central $\bbZ / 2$ subgroup are clear. Since we have a central extension
\[
1 \to \langle z_0 \rangle \cong \bbZ / 2 \to \SL_2(\calO) \to \PSL_2(\calO) \to 1,
\]
lifting each relation $\wh{R}_j$ from Theorem~\ref{thm:PresentPSL} to a word in $\sig$, $\mu$, $\tau$, and $\eta$ must give either the identity or $z_0$. This determines the given relation $R_j$.

We claim that $\SL_2(\calO)$ has no other relations. Indeed, the abstract group $\Lam$ with presentation in the statement of the corollary defines a central extension of $\PSL_2(\calO)$ by $\bbZ / 2$. Since these relations also hold in $\SL_2(\calO)$, the projection $\Lam \to \PSL_2(\calO)$ defined by killing $z_0$ factors through a homomorphism onto $\SL_2(\calO)$ defined by setting the abstract generators equal to the given matrix. Since this kernel has order two, $\SL_2(\calO)$ is isomorphic to one of $\Lam$ or $\PSL_2(\calO)$, and it is definitely not isomorphic to the latter.
\end{pf}
%%%%%%%%%%%%%%%%%%%%

%%%%%%%%%%%%%%%%%%%%
\begin{rem}
Observe that $\SL_2(\calO) \to \PSL_2(\calO)$ does not split. One way to see this is using the well-known fact that $\SL_2(\bbZ) \to \PSL_2(\bbZ)$ does not split.
\end{rem}
%%%%%%%%%%%%%%%%%%%%

Let $\Gam$ be the monodromy group of the Wiman--Edge pencil. We identify $\Gam$ with its image under the monodromy representation $\rho$ in \cite{FL1}, which is a subgroup of $\SL_2(\calO_o)$. Farb and Looijenga gave the following matrix generators.

%%%%%%%%%%%%%%%%%%%%
\begin{lem}\label{lem:WEwords}
The Wiman--Edge monodromy $\Gam$ is generated by:
\begin{align*}
\gam_\al &= \begin{pmatrix} 1 & -1 + 2 X \\ 0 & 1 \end{pmatrix} & \gam_{\al^\prime} &= \begin{pmatrix} 1 & 0 \\ 1 - 2 X & 1 \end{pmatrix} \\
&= \tau^{-1} \eta^2 & &= \sig \tau^{-1} \eta^2 \sig^{-1} \\
& & & \\
\gam_\beta &= \begin{pmatrix} 1 + X^3 & X^3 \\ - X^3 & 1 - X^3 \end{pmatrix} & \gam_{\beta^\prime} &= \begin{pmatrix} 1 + X^{-3} & X^{-3} \\ - X^{-3} & 1 - X^{-3} \end{pmatrix} \\
&= \tau^2 \eta^{-2} \sig \mu^3 \eta^{-2} \tau^4 & &= \eta^{-2} \tau^{-2} \sig \mu^{-3} \eta^{-2}
\end{align*}
\end{lem}
%%%%%%%%%%%%%%%%%%%%

%%%%%%%%%%%%%%%%%%%%
\begin{pf}
See \cite[Cor.~4.3]{FL1} for the matrix representatives for the generators. One then verifies by hand that the given words in our generators for $\SL_2(\calO)$ multiply out to the appropriate matrix.
\end{pf}
%%%%%%%%%%%%%%%%%%%%

%%%%%%%%%%%%%%%%%%%%
\begin{rem}
One can obtain relations for $\Gam$ using a computer algebra program like Magma \cite{Magma}. We found the relations to be sufficiently complicated that their inclusion wouldn't add any value to one's understanding of either $\Gam$ or the geometry of any of the objects associated with it.
\end{rem}
%%%%%%%%%%%%%%%%%%%%

%%%%%%%%%%%%%%%%%%%%
\section{Precise determination of the monodromy}\label{sec:Mono}
%%%%%%%%%%%%%%%%%%%%

The goal of this section is to give a more refined version of the following theorem of Farb--Looijenga.

%%%%%%%%%%%%%%%%%%%%
\begin{thm}[Thm.~1.1 \cite{FL1}]\label{thm:Arith}
The monodromy group of the Wiman--Edge pencil is isomorphic to a finite index subgroup of $\SL_2(\calO_o)$; in particular it is arithmetic.
\end{thm}
%%%%%%%%%%%%%%%%%%%%

Our proof will be independent of \cite{FL1} beyond our used of their matrix generators in Lemma~\ref{lem:WEwords}. Farb and Looijenga also proved that the reduction of $\Gam$ modulo $\frakp_5$ maps onto a unipotent subgroup of $\SL_2(\bbF_5)$ \cite[Cor.~4.3]{FL1}. They asked whether this congruence identity completely determines $\Gam$ \cite[Qn.~4.4]{FL1}, and our results answer this question in the negative. In particular, we will see that there is an additional congruence condition modulo $4 \calO$.

Before stating our precise determination of the monodromy, we need a preliminary lemma that describes $\SL_2(\calO / 4 \calO)$ as a short exact sequence of algebraic groups over $\bbF_4$. Observe first that $\calR_4 = \calO / 4 \calO$ can be described as $(\bbZ / 4)[X]$ and reduction modulo $2 \calR_4$ is $\bbF_2[X] \cong \bbF_4$.

%%%%%%%%%%%%%%%%%%%%
\begin{lem}\label{lem:Mod4}
There is a natural short exact sequence
\[
1 \to \fraksl_2(\bbF_4) \to \SL_2(\calR_4) \to \SL_2(\bbF_4) \to 1
\]
where $\fraksl_2(\bbF_4)$ denotes the Lie algebra of $\SL_2(\bbF_4)$, considered as an additive group. The image $\SL_2(\calO_o / 4 \calO)$ of $\SL_2(\calO_o)$ in $\SL_2(\calR_4)$ is the pullback of $\SL_2(\bbF_2) < \SL_2(\bbF_4)$ under this sequence.
\end{lem}
%%%%%%%%%%%%%%%%%%%%

%%%%%%%%%%%%%%%%%%%%
\begin{pf}
The first statement is very standard and follows from studying matrices in $\SL_2(\calR_4)$ of the form $\Id + 2 M$ and noticing that these are naturally represented by $M \in \fraksl_2(\bbF_4)$, where we consider $\bbF_4$ as $(\bbZ / 2)[X]$ inside $\calR_4 \cong (\bbZ/4)[X]$. The second statement follows from the fact that $\calO_o$ is the pullback of $\bbF_2 \subset \bbF_4$. 
\end{pf}
%%%%%%%%%%%%%%%%%%%%

%%%%%%%%%%%%%%%%%%%%
\begin{rem}
Since $\bbF_4$ has characteristic $2$, it can be identified with the additive group of matrices in $\M_2(\bbF_4)$ of the form $\left(\begin{smallmatrix} x_1 & x_2 \\ x_3 & x_1 \end{smallmatrix}\right)$, i.e., matrices with trace zero. As a group, $\fraksl_2(\bbF_4) \cong (\bbZ / 2)^6$.
\end{rem}
%%%%%%%%%%%%%%%%%%%%

We can now restate our main result.

%%%%%%%%%%%%%%%%%%%%
\begin{thma}\label{thm:Description}
The monodromy group $\Gam$ of the Wiman--Edge pencil can described equivalently in $\SL_2(\calO)$ and $\SL_2(\calO_o)$ by the following two properties:
\begin{enumerate}

\item The reduction of $\Gam$ modulo $\frakp_5$ is unipotent.

\item The reduction $\conj{\Gam}$ of $\Gam$ modulo $4 \calO$ fits into an exact sequence
\[
1 \to C_4 \to \conj{\Gam} \to \SL_2(\bbF_2) \to 1
\]
where $C_4 < \fraksl_2(\bbF_4)$ denotes the subgroup of matrices $\left(\begin{smallmatrix} x_1 & x_2 \\ x_3 & x_1 \end{smallmatrix}\right) \in \fraksl_2(\bbF_4)$ so that $\Tr_{\bbF_4 / \bbF_2}(x_1 + x_2 + x_3) = 0$, where $\Tr_{\bbF_4 / \bbF_2}$ denotes the trace.

\end{enumerate}
\end{thma}
%%%%%%%%%%%%%%%%%%%%

%%%%%%%%%%%%%%%%%%%%
\begin{pf}
The proof is carried out by interpreting computations done in the computer algebra system Magma \cite{Magma}. In particular, we describe how one gives arithmetic significance to standard Magma calculations and refer the reader to \cite{code} for code that can be used to execute these calculations. Using the presentation for $\SL_2(\calO)$ in Corollary~\ref{cor:PresentSL} and the generators for $\Gam$ from Lemma~\ref{lem:WEwords}, one obtains a presentation of $\Gam$ as a finite index subgroup of $\SL_2(\calO)$ with index $480$.

The fact that $\Gam$ has unipotent image modulo $\frakp_5$ is contained in \cite[Cor.~4.3]{FL1}, but we give an alternate argument as a warm-up for taking reductions modulo $2$ and $4$. One can check with Magma that there is a unique normal subgroup of $\SL_2(\calO)$ with quotient isomorphic to $\SL_2(\bbF_5)$. Therefore this subgroup must be the congruence subgroup $\SL_2(\calO)[\frakp_5]$ consisting of those matrices congruent to the identity modulo $\frakp_5$. The image of $\Gam$ in the quotient group has order $5$, which confirms our assertion, since all subgroups of $\SL_2(\bbF_5)$ isomorphic to $\bbZ / 5$ are unipotent.

Similarly, one sees that there are two normal subgroups of $\Gam$ with quotient isomorphic to $\SL_2(\bbF_4)$, where the second arises from its exceptional isomorphism with $\PSL_2(\bbF_5)$. The subgroup that does not contain $\SL_2(\calO)[\frakp_5]$ must be the congruence subgroup $\SL_2(\calO)[2]$ of level $2 \calO$. One checks that the image of $\Gam$ in the quotient is isomorphic to $\PSL_2(\bbF_2)$. Since $\Gam \le \SL_2(\calO_o)$, this image must be the standard $\SL_2(\bbF_2)$ in $\SL_2(\bbF_4)$.

We now study the reduction modulo $4 \calO$. From Lemma~\ref{lem:Mod4}, one sees that $\SL_2(\calR_4)$ has order $3840$. Magma confirms that there is a unique normal subgroup of $\SL_2(\calO)$ with that index, hence it must be the level $4$ congruence subgroup $\SL_2(\calO)[4]$. Note that $\fraksl_2(\bbF_4)$ from Lemma~\ref{lem:Mod4} is the image of $\SL_2(\calO)[2]$ under reduction modulo $4 \calO$.

One can then compute the intersection of $\fraksl_2(\bbF_4)$ with the image $\conj{\Gam}$ of $\Gam$, and one obtains an index two subgroup of $\fraksl_2(\bbF_4)$. Note that the group $C_4$ from the statement of the theorem is also an index two subgroup, since it is defined by a linear equation over $\bbZ / 2$, again considering $\bbF_4$ as $(\bbZ / 2)[X]$. One does not need Magma to see that $\mu^3$ reduces modulo $4 \calO$ to an element of $\fraksl_2(\bbF_4)$ not in $C_4$, hence it represents the nontrivial coset representative for $C_4$ in $\fraksl_2(\bbF_4)$. One then checks in Magma that the image of $\mu^3$ in $\fraksl_2(\bbF_4)$ is also not in $\fraksl_2(\bbF_4) \cap \conj{\Gam}$, and it follows that $\fraksl_2(\bbF_4) \cap \conj{\Gam} = C_4$, as desired.

The above verifies that $\Gam$ satisfies the congruence conditions given in the statement of the theorem. We must check that the above conditions completely describe $\Gam$ as a subgroup of $\SL_2(\calO)$. This is equivalent to showing that $\Gam$ contains the congruence subgroup $\SL_2(\calO)[4 \frakp_5]$. Magma can easily check this by showing that $\Gam$ contains the intersection of $\SL_2(\calO)[\frakp_5]$ and $\SL_2(\calO)[4]$. This completes our description of how one verifies the theorem in Magma.
\end{pf}
%%%%%%%%%%%%%%%%%%%%

%%%%%%%%%%%%%%%%%%%%
\section{The action on $\frakh \times \frakh$}\label{sec:Quo}
%%%%%%%%%%%%%%%%%%%%

Let $\frakh$ be the upper-half plane. We now study the action of $\Gam$ on $\frakh \times \frakh$ and the geometry of the quotient space $X_\Gam = (\frakh \times \frakh) / \Gam$. Note that it is $\PSL_2(\calO)$ that acts faithfully on $\frakh \times \frakh$, so it will be convenient to replace $\Gam$ with its image in $\PSL_2(\calO)$.

We first notice that $X_\Gam$ is a manifold, by which we mean as a Riemannian orbifold modeled on $\frakh \times \frakh$, not in the weaker sense of the analytic quasi-projective variety underlying $X_\Gam$ being nonsingular. This is well-known to follow immediately from the following.

%%%%%%%%%%%%%%%%%%%%
\begin{lem}\label{lem:TF}
The monodromy group $\Gam$ is torsion-free.
\end{lem}
%%%%%%%%%%%%%%%%%%%%

%%%%%%%%%%%%%%%%%%%%
\begin{pf}
Suppose $\gam \in \SL_2(\calO)$ has finite order. Then $\gam$ has characteristic polynomial $t^2 - \tr(\gam) t + 1$ and some root of unity $\zeta$ is a root of this polynomial. This implies that $k(\zeta)$ is at most a quadratic extension of $k$. This reduces us to
\[
n \in \{2, 3, 4, 5, 6, 10\}.
\]
An element of order $4$ must have square $z_0$, and $z_0 \notin \Gam$, so we can eliminate that case.

Now note that the characteristic polynomial $p_n(t)$ over $k$ for an element of order $n$ is
\[
p_n(t) = \begin{cases}
(t+1)^2 & n = 2 \\
t^2 + t + 1 & n = 3 \\
t^2 + X t + 1 & n = 5 \\
t^2 - t + 1 & n = 6 \\
t^2 - X t + 1 & n = 10
\end{cases}
\]
(For $n = 5$, and $n = 10$ there is technically also the Galois conjugate polynomial $t^2 \pm (1 - X) t + 1$.) Since the reduction of $\Gam$ modulo $1 - 2 X$ is unipotent, we must have that the reduction of $p_n(t)$ modulo $1 - 2 X$ factors as $(t-1)^2$. This reduces us immediately to the case $n = 5$.

Then, $p_5(t)$ does not factor over $\bbF_4 = \bbF_2(X)$. This implies that the congruence subgroup of level $2$ in $\SL_2(\calO)$ contains no $5$-torsion. Thus if $\Gam$ were to contain an element of order five, then its reduction modulo $2$ would be a nontrivial element of $\SL_2(\bbF_4)$ of order five. However, the reduction of $\Gam$ is $\SL_2(\bbF_2)$ has order six, so it contains no element of order five. This proves that $\Gam$ is torsion-free.
\end{pf}
%%%%%%%%%%%%%%%%%%%%

%%%%%%%%%%%%%%%%%%%%
\begin{rem}
Eduard Looijenga also described a nice independent proof of Lemma~\ref{lem:TF} using the geometry of the Wiman--Edge pencil. Briefly, the Torelli theorem implies that torsion in $\Gam$ would define a curve in the pencil with an exceptional automorphism centralizing the action of $A_5$, but there is no such curve.
\end{rem}
%%%%%%%%%%%%%%%%%%%%

We also have the following standard calculation.

%%%%%%%%%%%%%%%%%%%%
\begin{lem}\label{lem:OpenEuler}
The space $X_\Gam$ has topological Euler number $e(X_\Gam) = 16$.
\end{lem}
%%%%%%%%%%%%%%%%%%%%

%%%%%%%%%%%%%%%%%%%%
\begin{pf}
Let $\conj{\Gam}$ be the image of $\Gam$ in $\PSL_2(\calO)$. Since $\PSL_2(\calO)$ acts faithfully on $\frakh \times \frakh$, we have
\[
e(X_\Gam) = [\PSL_2(\calO) : \conj{\Gam}]\, e((\frakh \times \frakh) / \PSL_2(\calO)),
\]
where $e((\frakh \times \frakh) / \PSL_2(\calO))$ denotes the Euler--Poincar\'e characteristic. However, it is a classical fact that
\[
e((\frakh \times \frakh) / \PSL_2(\calO)) = 2 \zeta_k(-1) = \frac{1}{15},
\]
where $\zeta_k$ denotes the Dedekind zeta function; e.g., see \cite[Thm.~IV.1.1]{vdG}. Then $\Gam$ has index $480$ in $\SL_2(\calO)$ and intersects the center of $\SL_2(\calO)$ trivially (e.g., by Lemma~\ref{lem:TF}), hence $[\PSL_2(\calO) : \conj{\Gam}] = 240$ and the lemma follows.
\end{pf}
%%%%%%%%%%%%%%%%%%%%

Let $\Del < \SL_2(\calO)$ be the subgroup of upper-triangular matrices and $\conj{\Del}$ denote the image of $\Del$ in $\PSL_2(\calO)$. Since $k$ has class number one, it follows that $(\frakh \times \frakh) / \SL_2(\calO)$ has one cusp, i.e., $\Del$ determines the unique conjugacy class of parabolic subgroups of $\SL_2(\calO)$ \cite[Prop.~I.1.1]{vdG}.

Before counting the cusps of $X_\Gam$, we prove an easy lemma on the structure of $\Del$. The result is very well-known in much greater generality (e.g., see \cite[\S II.1]{vdG}), but we will need some calculations related to the proof later so we include it for the reader's convenience.

%%%%%%%%%%%%%%%%%%%%
\begin{lem}\label{lem:SLcusp}
There is a natural isomorphism $\Del \cong \conj{\Del} \times \bbZ / 2$, where the $\bbZ / 2$ factor is generated by the center $\langle z_0 \rangle$ of $\SL_2(\calO)$. The group $\conj{\Del}$ is generated by $\wh{\mu}, \wh{\tau}, \wh{\eta} \in \PSL_2(\calO)$ and has presentation on these generators with relations:
\begin{align*}
\wh{S}_0 &= [\wh{\tau}, \wh{\eta}] & & \\
\wh{S}_1 &= \wh{\mu} \wh{\tau} \wh{\mu}^{-1} (\wh{\tau} \wh{\eta})^{-1} & \wh{S}_2 &= \wh{\mu} \wh{\eta} \wh{\mu}^{-1} (\wh{\tau} \wh{\eta}^2)^{-1}
\end{align*}
\end{lem}
%%%%%%%%%%%%%%%%%%%%

%%%%%%%%%%%%%%%%%%%%
\begin{pf}
We start by presenting $\conj{\Del}$. Notice that $\wh{S}_0, \wh{S}_1, \wh{S}_2$ are just the relations $\wh{R}_4, \wh{R}_5, \wh{R}_6$ from Theorem~\ref{thm:PresentPSL}. It is also clear from the semi-direct product structure on the upper-triangular subgroup of $\PSL_2(\bbR)$ that we obtain a split short exact sequence:
\begin{equation}\label{eq:CuspSeq}
1 \to \langle \wh{\tau}, \wh{\eta} \rangle \cong \bbZ^2 \to \conj{\Del} \cong \langle \wh{\tau}, \wh{\eta} \rangle \rtimes \langle \mu \rangle \to \bbZ \to 1
\end{equation}
In particular, one must only understand the conjugation action of $\wh{\mu}$ on $\wh{\tau}$ and $\wh{\eta}$ to obtain a presentation for $\wh{\Del}$. This is precisely what is described by the relations $\wh{S}_1$ and $\wh{S}_2$. Finally, one checks that lifting to $\mu, \tau, \sig \in \Del$ splits the central exact sequence
\[
1 \to \langle z_0 \rangle \cong \bbZ / 2 \to \Del \to \conj{\Del} \to 1
\]
which proves that $\Del \cong \conj{\Del} \times \bbZ / 2$.
\end{pf}
%%%%%%%%%%%%%%%%%%%%

We now determine the cusps of $X_\Gam$ and identify the associated cusp subgroups of $\Gam$ as subgroups of $\Del$.

%%%%%%%%%%%%%%%%%%%%
\begin{prop}\label{prop:Cusps}
The space $X_\Gam$ has $6$ cusps. The associated conjugacy classes $\Del_1, \dots, \Del_6$ of parabolic subgroups of $\Gam$ are conjugate in $\SL_2(\calO)$ to the following subgroups of $\Del$:
\begin{align*}
\Del_1, \Del_2 \sim \Lam_8 &= \langle \mu^{-2} \eta, \eta^2, \tau^2 \rangle \\
\Del_3 \sim \Lam_{24} &= \langle \mu^6, \tau^2, \tau \eta^2 \rangle \\
\Del_4, \Del_5 \sim \Lam_{40} &= \langle \mu^2 \tau^{-2} \eta^{-1}, \tau^2 \eta^6, \eta^{10} \rangle \\
\Del_6 \sim \Lam_{120} &= \langle \mu^6, \tau \eta^{18}, \eta^{20} \rangle
\end{align*}
Here $\sim$ denotes conjugacy in $\SL_2(\calO)$ and $[\Del : \Lam_j] = j$. Specifically, we have:
\begin{align*}
\Del_1 &= (\sig \tau \mu^2)^{-1} \Lam_8 (\sig \tau \mu^2) &
\Del_2 &= (\tau \sig \tau \sig \mu^{-1}) \Lam_8 (\tau \sig \tau \sig \mu^{-1}) \\
\Del_3 &= (\sig \tau)^{-1} \Lam_{24} (\sig \tau) &
\Del_4 &= (\sig \eta)^{-1} \Lam_{40} (\sig \eta) \\
\Del_5 &= (\mu \sig \tau \mu)^{-1} \Lam_{40} (\mu \sig \tau \mu) &
\Del_6 &= \Lam_{120}
\end{align*}
\end{prop}
%%%%%%%%%%%%%%%%%%%%

%%%%%%%%%%%%%%%%%%%%
\begin{pf}
Again, we explain how one deduces the result using Magma. At the risk of some confusion, we work in $\PSL_2(\calO)$ instead of $\SL_2(\calO)$ but discard the $\wh{g}$ notation in favor of just $g$ for readability. For an ideal $\calJ$ of $\calO$, let
\[
X[\calJ] = (\frakh \times \frakh) / \PSL_2(\calO)[\calJ]
\]
denote the quotient of $\frakh \times \frakh$ by the principal congruence subgroup of $\PSL_2(\calO)$ of level $\calJ$.

It is an easy consequence of the orbit-stabilizer theorem that $X[\calJ]$ has $[\PSL_2(\calO / \calJ) : r_\calJ(\Del)]$ cusps, where
\[
r_\calJ : \PSL_2(\calO) \to \PSL_2(\calO / \calJ)
\]
is the reduction homomorphism. For $\calJ = 4 \frakp_5$ this index is $240$. One can then count cusps by identifying them with right coset representatives of $r_\calJ(\Del)$ in $\PSL_2(\calO / \calJ)$. Then $\Gam$ acts on these cosets through $r_\calJ$, and we can identify cusps of $X_\Gam$ with $\Gam$-orbits under this action.

Specifically, if we write
\[
\PSL_2(\calO / \calJ) = \coprod_{j = 1}^{240} r_\calJ(\Del) \ep_j,
\]
let $[\ep_j]$ be the cusp associated with $\ep_j$. Then $h \in r_\calJ(\Gam)$ acts on the right by $[\ep_j] \cdot h = [\ep_{h(j)}]$, where $\ep_j h = \del(h, j) \ep_{h(j)}$ for $\del(h, j) \in r_\calJ(\Del)$. The stabilizer in $\PSL_2(\calO / \calJ)$ of $[\ep_j]$ is $\ep_j^{-1} r_\calJ(\Del) \ep_j$, and it follows that $\ep_{j_1}$ is in the same $r_\calJ(\Gam)$-orbit as $\ep_{j_2}$ if and only if there is some $h \in r_\calJ(\Gam)$ so that $\ep_{j_1} h \ep_{j_2}^{-1} \in r_\calJ(\Del)$. A simple for loop finds that there are six orbits, hence $X_\Gam$ has six cusps.

Using an abstract presentation for $\Del$ derived from taking the quotient of the presentation in Lemma~\ref{lem:SLcusp} by the central $\bbZ / 2$, one finds subgroups $\Lam_{k(j)}$ of $\Del$ for which some conjugate of $\Lam_{k(j)}$ in $\Del$ is the pullback to $\Del$ of
\begin{equation}\label{eq:CuspMess}
\ep_j\left(\ep_j^{-1} r_\calJ(\Del) \ep_j \cap r_\calJ(\Gam) \right) \ep_j^{-1}.
\end{equation}
The groups in Equation~\eqref{eq:CuspMess} are the images under $r_\calJ$ of $\Lam_8$, $\tau^{-1} \Lam_8 \tau$, $\Lam_{24}$, $\Lam_{40}$, $\mu^{-1} \Lam_{40} \mu$, and $\Lam_{120}$. Taking a representative $e_j \in \SL_2(\calO)$ for $\ep_j$, we see that $ (e_j g)^{-1} \Lam_{k(j)} (e_j g)$ is a subgroup of $\Gam$ representing the cusp associated with $\ep_j$. Finding generators for each $\Lam_j$ in a standard way completes the proof.
\end{pf}
%%%%%%%%%%%%%%%%%%%%

%%%%%%%%%%%%%%%%%%%%
\section{The cusp resolutions}\label{sec:Resolve}
%%%%%%%%%%%%%%%%%%%%

We retain the notation of \S\ref{sec:Quo}, and now compute the standard smooth compactification $Y_\Gam$ of $X_\Gam$. We refer the reader to \cite[Ch.~II]{vdG} for more details of the construction of the compactification. We briefly sketch some general facts before moving on to describing each cusp resolution for $X_\Gam$ in detail. We work in $\PSL_2(\calO)$, but discard the $\,\widehat{\, }\,$ notation for readability. At the conclusion of this section, we will prove:

%%%%%%%%%%%%%%%%%%%%
\begin{thmb}\label{thm:ChernNos}
The smooth compactification $Y_\Gam$ of $X_\Gam$ is a smooth projective surface of general type with Chern numbers:
\begin{align*}
c_1^2(Y_\Gam) &= 16 \\
c_2(\Gam) &= 56
\end{align*}
It has holomorphic Euler characteristic $\chi(\calO_{Y_\Gam}) = 6$, irregularity $q(Y_\Gam) = 0$, and geometric genus $p_g(Y_\Gam) = 5$.
\end{thmb}
%%%%%%%%%%%%%%%%%%%%

Recall that we have the subgroup $\Del = \langle \mu, \tau, \eta \rangle$ of $\PSL_2(\calO)$. If $\Lam \le \Del$ is a finite index subgroup, then we can use the exact sequence in Equation~\eqref{eq:CuspSeq} to see that
\[
\Lam = \left\langle \mu^a y\,,\, \tau^{b_1} \eta^{b_2}\,,\, \tau^{c_1} \eta^{c_2} \right\rangle
\]
for $a, b_1, b_2, c_1, c_2 \in \bbZ$, some $y \in \langle \tau, \eta \rangle$, and $a \neq 0$. Moreover, we can choose these generators so that $\Lam$ has index $|a (b_1 c_2 - b_2 c_1)|$ in $\Del$. This implies that
\[
T_\Lam = \langle \tau^{b_1} \eta^{b_2}\,,\, \tau^{c_1} \eta^{c_2} \rangle = \Lam \cap \langle \tau, \eta\rangle
\]
is the kernel of the map from $\Lam$ to $\bbZ$ induced by Equation~\eqref{eq:CuspSeq}. Note that the generators for each $\Lam_j$ in the statement of Proposition~\ref{prop:Cusps} are chosen to satisfy these properties.

We now describe how to resolve a cusp associated with $\Lam$ as above. Consider the standard embedding of $k$ into $\bbR^2$ using its two archimedean places, which induces a discrete embedding of $T_\Lam \le \bbZ[X]$ into $\bbR^2$. Let $T_\Lam^+$ be the elements of $T_\Lam$ lying in the first quadrant of $\bbR^2$ and $\calC_\Lam$ be the convex hull of $T_\Lam$.

\begin{figure}
\centering
\begin{tikzpicture}[scale=0.7]
\draw (0,0) -- (10,0);
\draw (0,0) -- (0,10);
\draw[fill=blue!10] (10, 10) -- (3,3) -- (5,2) -- (10,4);
\draw[fill=green!10] (10,10*1.5/7) -- (7, 1.5) -- (5,2) -- (10,4);
\draw[fill=blue!10]  (10,10*1.5/7) -- (7, 1.5) -- (9, 1.25) -- (10, 12.5/9);
\draw[fill=green!10] (10, 1.2) -- (9, 1.25) -- (10, 12.5/9);
\draw[fill=green!10] (10, 10) -- (3,3) -- (2,5) -- (4,10);
\draw[fill=blue!10] (10*1.5/7,10) -- (1.5,7) -- (2,5) -- (4,10);
\draw[fill=green!10] (10*1.5/7,10) -- (1.5,7) -- (1.25,9) -- (12.5/9,10);
\draw[fill=blue!10] (1.2,10) -- (1.25,9) -- (12.5/9,10);
\draw[dashed] (0,0) -- (3,3);
\draw[dashed] (0,0) -- (2,5);
\draw[dashed] (0,0) -- (1.5, 7);
\draw[dashed] (0,0) -- (1.25, 9);
\draw[dashed] (0,0) -- (5,2);
\draw[dashed] (0,0) -- (7,1.5);
\draw[dashed] (0,0) -- (9, 1.25);
\draw[black,fill=red] (3,3) circle (2pt);
\draw[black,fill=red] (5,2) circle (2pt);
\draw[black,fill=red] (2,5) circle (2pt);
\draw[black,fill=red] (7,1.5) circle (2pt);
\draw[black,fill=red] (1.5,7) circle (2pt);
\draw[black,fill=red] (9,1.25) circle (2pt);
\draw[black,fill=red] (1.25,9) circle (2pt);
\draw[red] [->] (1.6,7) [bend left]  to node[left] {$\mu^a z$} (7,1.6);
\draw[thick, ->] (10.1, 5) to (11, 5);
\draw[red] (12, 2.7) -- (12, 7.3);
\draw[red] (16, 2.7) -- (16, 7.3);
\draw[red] (11.7, 3) -- (16.3, 3);
\draw[red] (11.7, 7) -- (16.3, 7);
\draw[black,fill=black] (12,3) circle (1.5pt);
\draw[black,fill=black] (12,7) circle (1.5pt);
\draw[black,fill=black] (16,3) circle (1.5pt);
\draw[black,fill=black] (16,7) circle (1.5pt);
\end{tikzpicture}
\caption{Compactification by a square of rational curves}\label{fig:CuspCartoon}
\end{figure}
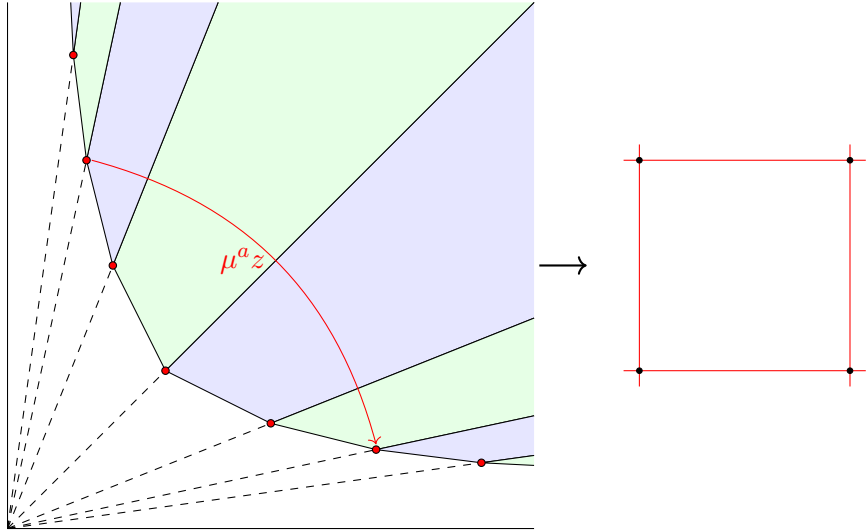

Then $\mu^a y$ acts on $T_\Lam$ by conjugation and preserves $T_\Lam^+$, hence it acts on $\calC_\Lam$. The action on the boundary $\calB_\Lam$ of $\calC_\Lam$ determines a combinatorial $n$-gon, where we associate edges of the $n$-gon with $(\mu^a y)$-orbits of vertices of $\calB_\Lam$ and vertices of the $n$-gon with edges of $\calB_\Lam$. See Figure~\ref{fig:CuspCartoon} for a general picture. Each edge of this $n$-gon represents a smooth rational curve in the compactification.

It remains to determine the self-intersection of the rational curve $E_v$ associated with a vertex $v$ on $\calB_\Lam$. The two adjacent vertices $v_1, v_{-1}$ of $\calB_\Lam$ have the property that, as elements of $\bbZ[X]$, $v_1 + v_{-1} = b v$ for some integer $b > 1$, and we have $E_v^2 = - b$. See \cite[\S II.2]{vdG}.

Finally, we note that how one resolves a cusp associated with some $\Lam \le \Del$ is independent of its conjugacy class in $\Del$. In particular, to resolve the cusps of $X_\Gam$, it suffices to study the groups $\Lam_j$ in Proposition~\ref{prop:Cusps}. The remainder of this section studies these resolutions, and Mathematica code to assist the reader in verifying assertions about various convex hulls is available on the author's webpage \cite{code2}.

%%%%%%%%%%%%%%%%%%%%
\subsection{Degree $8$ cusps}
%%%%%%%%%%%%%%%%%%%%

\begin{figure}
\centering
\begin{tikzpicture}
\draw (0,0) -- (10,0);
\draw (0,0) -- (0,10);
\draw[fill=blue!10] (10, 10) -- (3,3) -- (5,2) -- (10,4);
\draw[fill=green!10] (10,10*1.5/7) -- (7, 1.5) -- (5,2) -- (10,4);
\draw[fill=blue!10]  (10,10*1.5/7) -- (7, 1.5) -- (9, 1.25) -- (10, 12.5/9);
\draw[fill=green!10] (10, 1.2) -- (9, 1.25) -- (10, 12.5/9);
\draw[fill=green!10] (10, 10) -- (3,3) -- (2,5) -- (4,10);
\draw[fill=blue!10] (10*1.5/7,10) -- (1.5,7) -- (2,5) -- (4,10);
\draw[fill=green!10] (10*1.5/7,10) -- (1.5,7) -- (1.25,9) -- (12.5/9,10);
\draw[fill=blue!10] (1.2,10) -- (1.25,9) -- (12.5/9,10);
\draw[dashed] (0,0) -- (3,3);
\draw[dashed] (0,0) -- (2,5);
\draw[dashed] (0,0) -- (1.5, 7);
\draw[dashed] (0,0) -- (1.25, 9);
\draw[dashed] (0,0) -- (5,2);
\draw[dashed] (0,0) -- (7,1.5);
\draw[dashed] (0,0) -- (9, 1.25);
\draw[black,fill=red] (3,3) circle (2pt);
\draw[black,fill=red] (5,2) circle (2pt);
\draw[black,fill=red] (2,5) circle (2pt);
\draw[black,fill=red] (7,1.5) circle (2pt);
\draw[black,fill=red] (1.5,7) circle (2pt);
\draw[black,fill=red] (9,1.25) circle (2pt);
\draw[black,fill=red] (1.25,9) circle (2pt);
\node[red, below] at (3,3) {$2$};
\node[red, right, scale=0.9] at (2,5) {$4 - 2 X$};
\node[red, above, scale=0.9] at (5,2.2) {$2 + 2 X$};
\node[red, right, scale=0.8] at (1.5,7) {$10 - 6 X$};
\node[red, above, scale=0.8] at (7,1.6) {$4 + 6 X$};
\node[red, left, scale=0.7] at (1.25,9) {$26 - 16 X$};
\node[red, above, scale=0.7] at (9,1.35) {$10 + 16 X$};
\end{tikzpicture}
\caption{Resolving the index $8$ cusps by two $-3$ curves}\label{fig:Res8}
\end{figure}
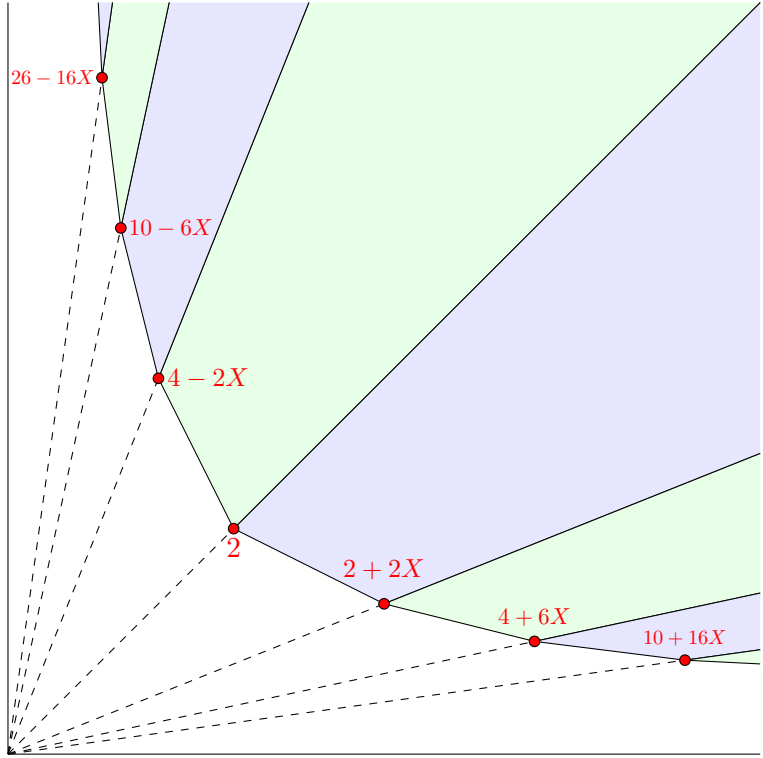

Here we study the resolution of a cusp associated with the group
\[
\Lam_8 = \langle \mu^{-2} \eta, \eta^2, \tau^2 \rangle.
\]
In the above notation, $T_\Lam = \langle \eta^2, \tau^2 \rangle = 2 \calO$ and $\mu^{-2} \eta$ has matrix
\[
\varphi = \begin{pmatrix} 5 & -3 \\ -3 & 2 \end{pmatrix}
\]
for its action on $T_\Lam$. The convex hull of $T_\Lam^+$ is pictured in Figure~\ref{fig:Res8}. Since $2$ is mapped to $4 - 6 X$ under $\varphi$, the cusp resolution is by a bigon of rational curves. Then
\begin{align*}
(10 - 6 X) + 2 &= 3(4 - 2 X) \\
(4 - 2 X) + (2 + 2 X) &= 3(2)
\end{align*}
which implies that both rational curves in the compactification have self-intersection $-3$.

%%%%%%%%%%%%%%%%%%%%
\subsection{Degree $24$ cusp}
%%%%%%%%%%%%%%%%%%%%

\begin{figure}[h]
\centering
\begin{tikzpicture}[scale=1.2]
\draw (0,10) -- (0,0) -- (10,0);
\draw[fill=blue!10] (10, 10) -- (0.5, 0.5) -- (1.55902, 0.440983) -- (10, 10*0.440983/1.55902);
\draw[fill=green!10] (10, 10) -- (0.5, 0.5) -- (0.440983,1.55902) -- (10*0.440983/1.55902,10);
\draw[fill=blue!20] (10, 10*0.440983/1.55902) -- (1.55902, 0.440983) -- (2.61803, 0.38196) -- (10, 10*0.38196/2.61803);
\draw[fill=green!20] (10*0.440983/1.55902,10) -- (0.440983,1.55902) -- (0.38196, 2.61803) -- (10*0.38196/2.61803, 10);
\draw[fill=blue!30] (10, 10*0.38196/2.61803) -- (2.61803, 0.38196) -- (3.67705, 0.32294) -- (10, 10*0.32294/3.67705);
\draw[fill=green!30] (10*0.38196/2.61803, 10) -- (0.38196, 2.61803)-- (0.32294, 3.67705) -- (10*0.32294/3.67705, 10);
\draw[fill=blue!40] (10, 10*0.32294/3.67705) -- (3.67705, 0.32294) -- (4.73606, 0.26393) -- (10, 10*0.26393/4.73606);
\draw[fill=green!40] (10*0.32294/3.67705, 10) -- (0.32294, 3.67705) -- (0.26393, 4.73606) -- (10*0.26393/4.73606, 10);
\draw[fill=blue!50] (10, 10*0.26393/4.73606) -- (4.73606, 0.26393) -- (5.79508, 0.20491) -- (10, 10*0.20491/5.79508);
\draw[fill=green!50] (10*0.26393/4.73606, 10) -- (0.26393, 4.73606) -- (0.20491, 5.79508) -- (10*0.20491/5.79508, 10);
\draw[fill=blue!60] (10, 10*0.20491/5.79508) -- (5.79508, 0.20491) -- (6.85410, 0.14589) -- (10, 10*0.14589/6.85410);
\draw[fill=green!60] (10*0.20491/5.79508, 10) -- (0.20491, 5.79508) -- (0.14589, 6.85410) -- (10*0.14589/6.85410, 10);
\draw[fill=blue!70] (10, 10*0.14589/6.85410) -- (6.85410, 0.14589) -- (7.91311, 0.08688) -- (10, 10*0.08688/7.91311);
\draw[fill=green!70] (10*0.14589/6.85410, 10) -- (0.14589, 6.85410) -- (0.08688, 7.91311) -- (10*0.08688/7.91311, 10);
\draw[fill=blue!80] (10, 10*0.08688/7.91311) -- (7.91311, 0.08688) -- (8.97213, 0.02786) -- (10, 10*0.02786/8.97213);
\draw[fill=green!80] (10*0.08688/7.91311, 10) -- (0.08688, 7.91311) -- (0.02786, 8.97213) -- (10*0.02786/8.97213, 10);
\foreach \Point in {(0.5, 0.5), (1.55902, 0.440983), (0.440983, 1.55902), (2.61803, 0.38196), (0.38196, 2.61803), (3.67705, 0.32294), (0.32294, 3.67705), (4.73606, 0.26393), (0.26393, 4.73606), (5.79508, 0.20491), (0.20491, 5.79508), (6.85410, 0.14589), (0.14589, 6.85410), (7.91311, 0.08688), (0.08688, 7.91311), (8.97213, 0.02786), (0.02786, 8.97213)
}{\draw[black,fill=red] \Point circle (1pt);};
\node[red, right, scale=0.75] at (0.4, 0.2) {$(2 + j) + 2 j X \quad j = 0,\dots,8$};
\node[red, right, scale=0.75] at (1.5, 5) {$(2 + 3 j) - 2 j X \quad j = 0,\dots,8$};
\end{tikzpicture}
\caption{Resolving the index $24$ cusp by sixteen curves}\label{fig:Res24}
\end{figure}
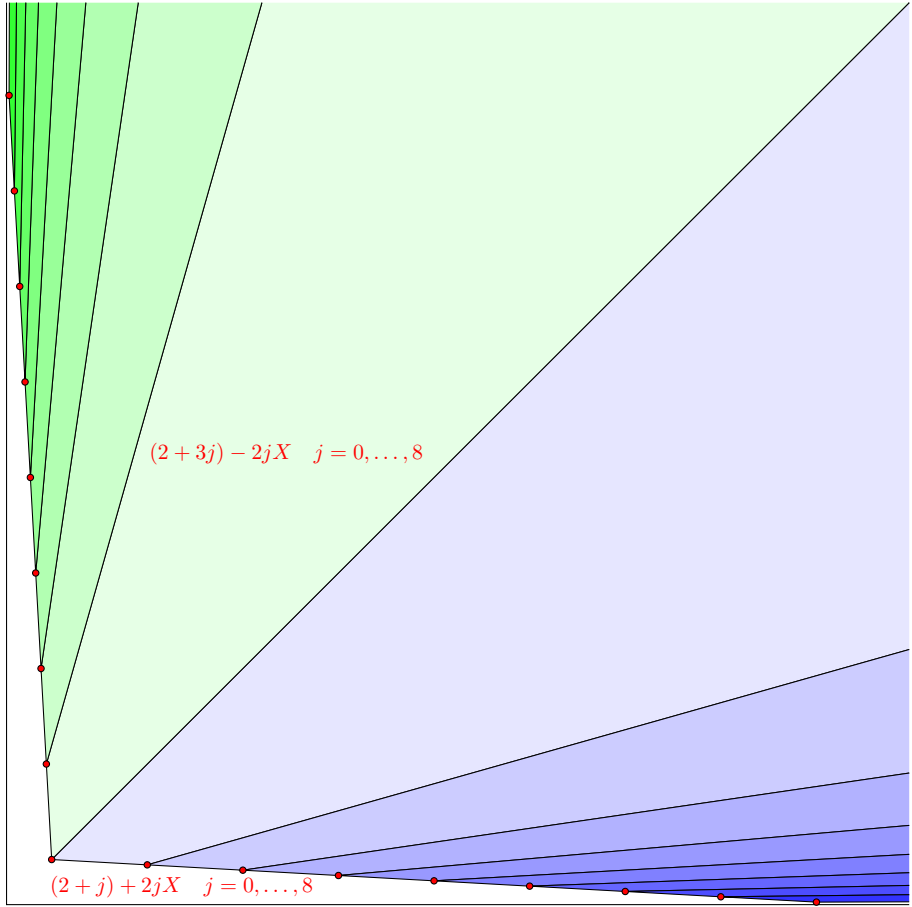

We now study the resolution of a cusp associated with the group
\[
\Lam_{24} = \langle \mu^6, \tau^2, \tau \eta^2 \rangle.
\]
In the above notation, $T_\Lam = \langle \tau^2, \tau \eta^2 \rangle$ and $\mu^6$ has matrix
\[
\varphi = \begin{pmatrix} 17 & 36 \\ 144 & 305 \end{pmatrix}
\]
for its action on $T_\Lam$. The convex hull of $T_\Lam^+$ is pictured in Figure~\ref{fig:Res24}. Then $\varphi$ sends $26 - 16 X$ to $10 + 16 X$, which implies that the compactification is by a $16$-gon of rational curves. For any adjacent vertices of the boundary of the convex hull of $T_\Lam^+$ that lie on a straight line, the associated curve has self-intersection $-2$, since:
\begin{align*}
(2 + j - 1 + 2 (j - 1) X) + (2 + j +1 + 2 (j + 1) X) &= 2 (2 + j + 2 j X) \\
(2 + 3(j - 1) - 2 (j - 1) X) + (2 + 3(j +1) - 2 (j + 1) X) &= 2 (2 + 3 j - 2 j X)
\end{align*}
There are fourteen such curves.

It remains to compute the self-intersections of the curves associated with $2$ and $10 + 16 X$. In these cases,
\begin{align*}
(3 + 2 X) + (5 - 2 X) &= 4(2)\\
(9 + 14X) + (31 + 50 X) & 4(10 + 16 X)
\end{align*}
which implies that the remaining two curves have self-intersection $-4$. In summary, the cusp associated with $\Lam_{24}$ is compactified by a $16$-gon of rational curves for which one pair of opposite sides determine curves with self-intersection $-4$ and the remaining curves have self-intersection $-2$.

%%%%%%%%%%%%%%%%%%%%
\subsection{Degree $40$ cusps}
%%%%%%%%%%%%%%%%%%%%

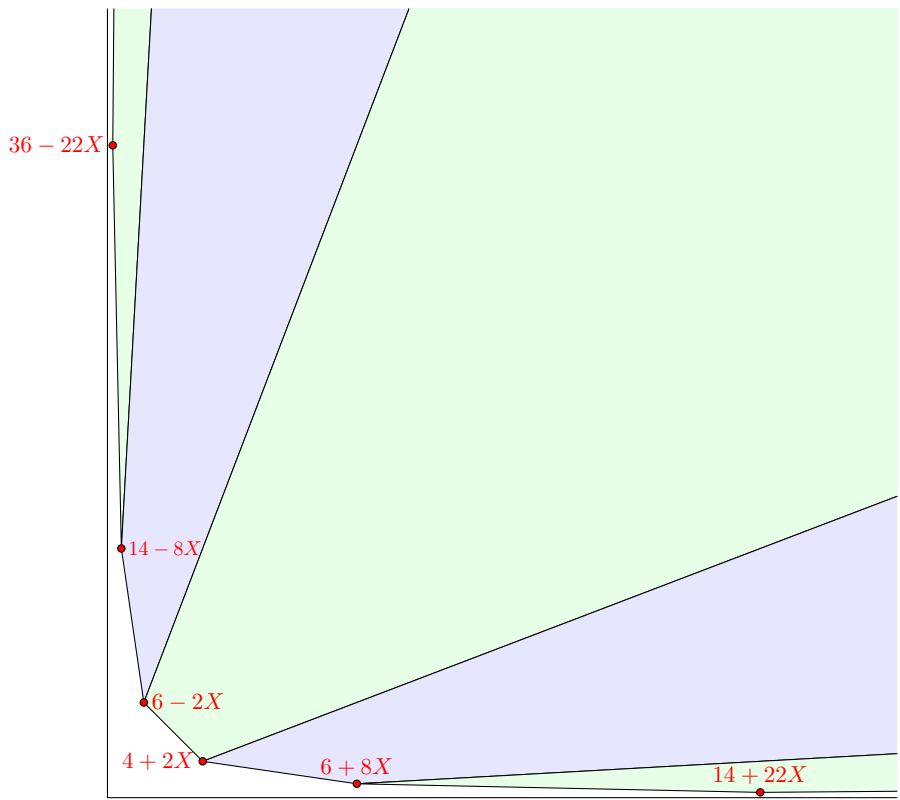
\begin{figure}[h]
\centering
\begin{tikzpicture}[scale=0.7]
\draw (0,15) -- (0,0) -- (15,0);
\draw[fill=green!10] (15*0.263932/4.73607, 15) -- (0.263932, 4.73607) -- (0.100813, 12.3992) -- (15*0.100813/12.3992, 15);
\draw[fill=blue!10] (15*0.690983/1.80902, 15) -- (0.690983, 1.80902) -- (0.263932, 4.73607) -- (15*0.263932/4.73607, 15);
\draw[fill=green!10] (15*0.690983/1.80902, 15) -- (0.690983, 1.80902) -- (1.80902, 0.690983) -- (15, 15*0.690983/1.80902);
\draw[green!10, fill=green!10] (15*0.690983/1.80902, 15) -- (15, 15) -- (15, 15*0.690983/1.80902);
\draw[fill=blue!10] (15, 15*0.690983/1.80902) -- (1.80902, 0.690983) -- (4.73607, 0.263932) -- (15, 15*0.263932/4.73607);
\draw[fill=green!10] (15, 15*0.263932/4.73607) -- (4.73607, 0.263932) -- (12.3992, 0.100813) -- (15, 15*0.100813/12.3992);
\foreach \Point in {
(12.3992, 0.100813), (4.73607, 0.263932), (1.80902, 0.690983), (0.690983, 1.80902), (0.263932, 4.73607), (0.100813, 12.3992)
}{\draw[black,fill=red] \Point circle (2pt);};
\node[red, above, scale=0.8] at (12.3992, 0.100813) {$14 + 22 X$};
\node[red, above, scale=0.8] at (4.73607, 0.263932) {$6 + 8 X$};
\node[red, left, scale=0.8] at (1.80902, 0.690983) {$4 + 2 X$};
\node[red, right, scale=0.8] at (0.690983, 1.80902) {$6 - 2 X$};
\node[red, right, scale=0.7] at (0.263932, 4.73607) {$14 - 8 X$};
\node[red, left, scale=0.8] at (0.100813, 12.3992) {$36 - 22 X$};
\end{tikzpicture}
\caption{Resolving the index $40$ cusps by two $-3$ curves}\label{fig:Res40}
\end{figure}

We assume the reader is becoming familiar with the nature of these calculations and start to skip some details. Considering the convex hull in Figure~\ref{fig:Res40}, for $\Lam_{40}$ we have that $\mu^2 \tau^{-2} \eta^{-1}$ sends $6 - 2 X$ to $6 + 8 X$. One then calculates that our cusp is compactified by a bigon of $-3$ curves.

%%%%%%%%%%%%%%%%%%%%
\subsection{Degree $120$ cusp}
%%%%%%%%%%%%%%%%%%%%

This case is completely analogous to the $\Lam_{24}$ case, and the cusp is compactified by a $16$-gon consisting of two $-4$ curves and fourteen $-2$ curves. Again, the $-4$ curves are on opposite sides of the $16$-gon. The boundary of the convex hull of $T_\Lam^+$ is given by the points:
\begin{align*}
(6 + j) + (8 - 2 j) X&& 0 \le j \le 8 \\
(6 + 11 j) + (8 + 18 j) X&& 0 \le j \le 8
\end{align*}
The lists overlap at $j = 0$, and the point $6 + 8 X$ represents one of the curves of self-intersection $-4$. At $j = 8$, we have $X^{12} (14 - 8 X) = 94 + 152 X$, which implies that $\mu^6$ identifies these two points; this represents the other curve of self-intersection $-4$. The picture of the convex hull is exactly as in Figure~\ref{fig:Res24}, appropriately relabeled.

%%%%%%%%%%%%%%%%%%%%
\subsection{Proof of Theorem~B}\label{ssec:ChernPf}
%%%%%%%%%%%%%%%%%%%%

The computation of the Chern numbers of $Y_\Gam$ follows from our above calculations of the cusp resolutions, Lemma~\ref{lem:OpenEuler}, and the formulas on \cite[p.63]{vdG}. It is well known that $q(Y_\Gam) = 0$ for any compactification of a quotient of $\frakh \times \frakh$ by an irreducible lattice (e.g., see \cite[Lem.~I.6.3]{vdG}). We then compute the remaining numerical invariants by:
\begin{align*}
\chi(\calO_{Y_\Gam}) &= \frac{1}{12}\left(c_1^2(Y_\Gam) + c_2(Y_\Gam)\right)\\
&= 1 - q(Y_\Gam) + p_g(Y_\Gam)
\end{align*}
Since $Y_\Gam$ is smooth and projective, it remains to show that $Y_\Gam$ is of general type. If it is minimal, this is clear from the possibilities for $c_1^2$ and $c_2$ in the Enriques--Kodaira classification \cite[Ch.~VI]{BPV}. However, performing a sequence of blowdowns on $Y_\Gam$ only increases $c_1^2$, and a minimal surface with $c_1^2 \ge 16$ must be of general type. This completes the proof. \qed

%%%%%%%%%%%%%%%%%%%%
\section{Image of the period map}\label{sec:Image}
%%%%%%%%%%%%%%%%%%%%

In this section, we study the image of the period map. Let $\calB^\circ$ be the base of the smooth locus of the pencil. Following Farb and Looijenga \cite[\S 5.1]{FL1}, we consider the period map $\Pi^\circ$ as a map from $\calB^\circ$ to
\[
X^\circ = (\frakh \times \frakh) / \SL_2(\calO_o)
\]
through the orbifold covering $X_\Gam \to X^\circ$. It is known that $\Pi^\circ$ is a closed embedding by \cite[Prop.~5.1]{FL1}.

\begin{figure}[h]
\centering
\begin{tikzcd}
& X^\circ \arrow[drr, hookleftarrow, "\Pi^\circ"] & & \\
& & X_\Gam \arrow[ul] \arrow[r, hookleftarrow] & \calB^\circ \\
X[2] \arrow[uur] & & & \\
& & X_\Gam[2] \arrow[uu] \arrow[ull] \arrow[r, hookleftarrow, dashed] & \wt{\calB}^\circ \arrow[uu] \\
\end{tikzcd}
\caption{Diagram of coverings for the period map} \label{fig:diagram}
\end{figure}
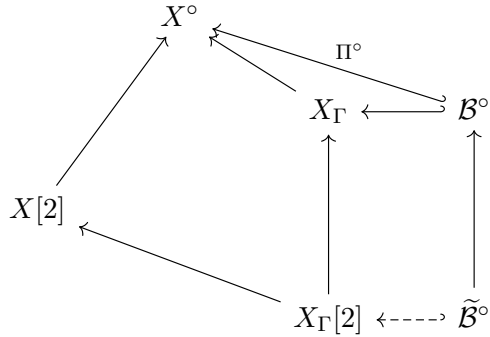

As described in \cite{FL1}, it is natural to consider $X^\circ$ as the quotient of
\[
X[2] = (\frakh \times \frakh) / \SL_2(\calO)[2]
\]
by an action of $\SL_2(\bbF_2)$, where $\SL_2(\calO)[2]$ denotes the level two congruence subgroup of $\SL_2(\calO)$. We use the classical connection between $X[2]$ and the Clebsch cubic surface to study the image of the period map. In other words, we will study the diagram of maps in Figure~\ref{fig:diagram}, where $X_\Gam[2]$ denotes the minimal common covering of $X_\Gam$ and $X[2]$ and $\wt{\calB}^\circ$ is the induced covering of $\calB^\circ$.

%%%%%%%%%%%%%%%%%%%%
\subsection{The curve on X[2]}\label{ssec:X2curve}
%%%%%%%%%%%%%%%%%%%%

Our first goal is to show that the image of the period map is associated with a (possibly singular) curve on $X[2]$ with genus one and $18$ punctures. We will then use the famous realization of a smooth compactification of $X[2]$ as a blowup of $\bbP^2$ to show that the image of the period map is determined by a plane curve of genus one that is invariant under an action of the symmetric group $S_3$. The following result immediately implies Theorem~C.

%%%%%%%%%%%%%%%%%%%%
\begin{prop}\label{prop:Bcover}
With notation as in Figure~\ref{fig:diagram}, $X_\Gam[2]$ is a Galois cover of $X_\Gam$ with group $\SL_2(\bbF_2)$. The cover $\wt{\calB}^\circ \to \calB^\circ$ is also Galois with group $\SL_2(\bbF_2)$, and $\wt{\calB}^\circ$ has genus one with $18$ punctures. The diagram of maps in Figure~\ref{fig:diagram} induces a closed embedding $\wt{\Pi}^\circ : \wt{\calB}^\circ \to X[2]$.
\end{prop}
%%%%%%%%%%%%%%%%%%%%

%%%%%%%%%%%%%%%%%%%%
\begin{pf}
The cover $X[2] \to X^\circ$ is Galois with group $\SL_2(\bbF_2)$ since $\SL_2(\calO_o)$ is the pull-back of $\SL_2(\bbF_2)$ to $\SL_2(\calO)$ under reduction modulo $2 \calO$. By Theorem~A, the monodromy $\Gam$ also has image $\SL_2(\bbF_2)$ under reduction modulo $2 \calO$. Elementary covering space theory then implies that the minimal common cover $X_\Gam[2]$ is a Galois cover of $X_\Gam$ with group $\SL_2(\bbF_2)$ that makes the left-most square in Figure~\ref{fig:diagram} into a commutative square of covers.

The map $\pi_1(\calB^\circ) \to \Gam$ is onto, hence the cover $\wt{\calB}^\circ$ of $\calB^\circ$ induced by $X_\Gam[2]$ is connected with group $\SL_2(\bbF_2)$. Since $\Pi^\circ$ is a closed embedding, we see that the induced map $\wt{\Pi}^\circ$ from $\wt{\calB}^\circ$ to $X[2]$ must also be a closed embedding. It remains to show that $\wt{\calB}^\circ$ is genus one with $18$ punctures.

Using the generators for $\pi_1(\calB^\circ)$ in Lemma~\ref{lem:WEwords} and defining reduction of $\Gam$ modulo $2 \calO$ as a homomorphism to $\SL_2(\bbF_2)$, the number of lifts of each puncture of $\calB^\circ$ to $\wt{\calB}^\circ$ is the index in $\SL_2(\bbF_2)$ of the subgroup generated by the image of an element associated with a loop around that puncture. The images of $\gam_\al$, $\gam_{\al^\prime}$, $\gam_\beta$, $\gam_{\beta^\prime}$ all have order $2$, so each associated puncture has three lifts. The element $\gam_\al \gam_\beta \gam_{\al^\prime} \gam_{\beta^\prime}$ is the identity modulo $2 \calO$, so the associated puncture has six lifts. We then have that $\wt{\calB}^\circ$ has Euler characteristic $-18$ and $18$ punctures, so it has genus one. This completes the proof.
\end{pf}
%%%%%%%%%%%%%%%%%%%%

%%%%%%%%%%%%%%%%%%%%
\subsection{The Clebsch cubic and the Klein plane}\label{ssec:Clebsch}
%%%%%%%%%%%%%%%%%%%%

We now describe the smooth compactification $Y[2]$ of $X[2]$ and its relationship to $\bbP^2$. See \cite[\S VIII.2]{vdG} for details. The space $Y[2]$ is the blowup of the famous Clebsch cubic surface at its ten Eckardt points, and the map
\[
\sig : Y[2] \to \bbP^2
\]
can be described as blowing up $\bbP^2$ at the vertices of an icosahedron and its dually inscribed dodecahedron. Farb and Looijenga call the image $\bbP^2$ the \emph{Klein plane}. In the remainder of this subsection we identify a number of the beautiful geometric features of this construction that we will use in studying the image of the period map.

%%%%%%%%%%%%%%%%%%%%
\subsubsection{Icosahedron vertices and Eckardt points}\label{sssec:Eckardt}
%%%%%%%%%%%%%%%%%%%%

Blowing up $\bbP^2$ at the six points associated with the vertices of the icosahedron gives the Clebsch surface. The ten points associated with vertices of the dodecahedron lift to the Eckardt points on the Clebsch surface. On $\bbP^2$, we denote the six icosahedral vertex points as $v_1, \dots, v_6$ and the dodecahedral vertex points associated with Eckardt points as $e_1, \dots, e_{10}$. See Figures~\ref{fig:IcosV} and~\ref{fig:DodecV} for coordinates.

\begin{figure}[h]
\begin{align*}
v_1 &= \begin{bmatrix} 0 \\ 1 \\ X \end{bmatrix} & v_3 &= \begin{bmatrix} X \\ 0 \\ 1 \end{bmatrix} & v_5 &= \begin{bmatrix} 1 \\ X \\ 0 \end{bmatrix} \\
v_2 &= \begin{bmatrix} 0 \\ 1 \\ -X \end{bmatrix} & v_4 &= \begin{bmatrix} -X \\ 0 \\ 1 \end{bmatrix} & v_6 &=\begin{bmatrix} 1 \\ -X \\ 0 \end{bmatrix}
\end{align*}
\caption{Vertices of the icosahedron}\label{fig:IcosV}
\end{figure}

\begin{figure}[h]
\begin{align*}
e_1 &= \begin{bmatrix} 1 \\ 1 \\ 1 \end{bmatrix} & e_3 &= \begin{bmatrix} 1 \\ -1 \\ 1 \end{bmatrix} & e_5 &= \begin{bmatrix} X-1 \\ 0 \\ X \end{bmatrix} & e_7 &= \begin{bmatrix} 0 \\ X \\ X-1 \end{bmatrix} & e_9 &= \begin{bmatrix} X \\ X-1 \\ 0 \end{bmatrix} \\
e_2 &= \begin{bmatrix} -1 \\ 1 \\ 1 \end{bmatrix} & e_4 &= \begin{bmatrix} -1 \\ -1 \\ 1 \end{bmatrix} &  e_6 &= \begin{bmatrix} 1-X \\ 0 \\ X \end{bmatrix} & e_8 &= \begin{bmatrix} 0 \\ X \\ 1-X \end{bmatrix} & e_{10} &= \begin{bmatrix} X \\ 1-X \\ 0 \end{bmatrix}
\end{align*}
\caption{Vertices of the dodecahedron}\label{fig:DodecV}
\end{figure}

As a Hilbert modular surface, the exceptional curves on $Y[2]$ arising from the Eckardt points are the compactifications of the ten lifts of $\frakh / \PSL_2(\bbZ)$ to $X[2]$. Notice that each lift is the thrice-punctured sphere arising from the quotient of $\frakh$ by the level $2$ congruence subgroup of $\PSL_2(\bbZ)$.

%%%%%%%%%%%%%%%%%%%%
\subsubsection{The compactification curves}\label{sssec:CuspCurves}
%%%%%%%%%%%%%%%%%%%%

Recall that the cusps of $X[2]$ are in one-to-one correspondence with elements of the projective like $\bbP^1(\bbF_4)$ over $\bbF_4$, or equivalently with the $\SL_2(\calO)[2]$-orbits of lifts of these lines to $\bbP^1(k)$. Each of the five cusps of $X[2]$ is compactified by a triangle of lines, where each line has self-intersection $-3$. Each line is the proper transform of the line of $\bbP^2$ associated with a pair of opposite edges of the icosahedron (or, dually, the dodecahedron). In particular, the image on $\bbP^2$ of each cusp line on $Y[2]$ meets exactly two of the six icosahedral vertices and exactly two of the dodecahedral vertices. We denote the line through $e_i$ and $e_j$ by $\ell_{i,j}$. See Figure~\ref{fig:DodecE} for an explicit parametrization of each line $\ell_{i,j}$ and Figure~\ref{fig:IcosE} for which $v_m$ lies on each $\ell_{i,j}$.

\begin{figure}[h]
\begin{align*}
\ell_{i, j} = \Bigg{\{} &r e_i + s e_j \ :\ \begin{bmatrix} r \\ s \end{bmatrix} \in \bbP^1 \Bigg{\}} \\
\{i,j\} \in \big{\{} &\{1,5\}\,,\, \{1,7\}\,,\, \{1,9\}\,,\, \{2,6\}\,,\, \{2,7\}\,,\, \{2,10\}\,,\, \\
&\{3,5\}\,,\, \{3,8\}\,,\, \{3,10\}\,,\, \{4,6\}\,,\, \{4,8\}\,,\, \{4,9\}\,,\, \\
&\{5,6\}\,,\, \{7,8\}\,,\, \{9,10\}\big{\}}
\end{align*}
\caption{Edges of the dodecahedron}\label{fig:DodecE}
\end{figure}

\begin{figure}[h]
\begin{align*}
v_2, v_5 &\in \ell_{1,5} & v_3, v_6 &\in \ell_{1,7} & v_1, v_4 &\in \ell_{1,9} \\
v_2, v_6 &\in \ell_{2,6} & v_4, v_5 &\in \ell_{2,7} & v_1, v_3 &\in \ell_{2,10} \\
v_1, v_6 &\in \ell_{3,5} & v_3, v_5 &\in \ell_{3,8} & v_2, v_4 &\in \ell_{3,10} \\
v_1, v_5 &\in \ell_{4,6} & v_4, v_6 &\in \ell_{4,8} & v_2, v_3 &\in \ell_{4,9} \\
v_3, v_4 &\in \ell_{5,6} & v_1, v_2 &\in \ell_{7,8} & v_5, v_6 &\in \ell_{9,10}
\end{align*}
\caption{Vertices of the icosahedron on the edge $\ell_{i,j}$}\label{fig:IcosE}
\end{figure}

%%%%%%%%%%%%%%%%%%%%
\subsubsection{The $\PSL_2(\bbF_2)$ action}\label{sssec:Action}
%%%%%%%%%%%%%%%%%%%%

The automorphism group of the Clebsch surface is the symmetric group $S_5$, and the induced action on $Y[2]$ restricts to the automorphism group of $X[2]$. This connects to the action of $\PSL_2(\bbF_4)$ on $X[2]$ through its isomorphism with the alternating group $A_5$. (The remaining order $2$ automorphism of $X[2]$ needed to generate $S_5$ is induced by the ``swap map'' exchanging the two factors of $\frakh \times \frak h$.) The inclusion of $\PSL_2(\bbF_2)$ into $\PSL_2(\bbF_4)$ is then realized as an inclusion of $S_3$ in $A_5$. See Figure~\ref{fig:S3} for matrices generating the $S_3$ action and Figures~\ref{fig:S3V} and~\ref{fig:S3E} for vertex and edge orbits under this action.

\begin{figure}[h]
\begin{align*}
\sig_1 &= \begin{bmatrix} -1 & 0 & 0 \\ 0 & -1 & 0 \\ 0 & 0 & 1 \end{bmatrix} & \sig_2 &= \begin{bmatrix} -\frac{X}{2} & -\frac{1}{2} & \frac{X-1}{2} \\[6pt] -\frac{1}{2} & \frac{X-1}{2} & -\frac{X}{2} \\[6pt] \frac{X-1}{2} & -\frac{X}{2} & -\frac{1}{2} \end{bmatrix} \\
\tau &= \begin{bmatrix} \frac{X}{2} & \frac{1}{2} & \frac{1-X}{2} \\[6pt] \frac{1}{2} & \frac{1-X}{2} & \frac{X}{2} \\[6pt] \frac{X-1}{2} & -\frac{X}{2} & -\frac{1}{2} \end{bmatrix} = \sig_1 \sig_2 & \sig_3 &= \tau \sig_1 = \sig_1 \sig_2 \sig_1
\end{align*}
\[
S_3 = \left\langle \sig_1, \sig_2, \tau\ |\ \sig_1^2, \sig_2^2, \tau^3, \tau=\sig_1 \sig_2 \right\rangle
\]
\caption{The $S_3$ action}\label{fig:S3}
\end{figure}

\begin{figure}[h]
\begin{align*}
\sig_1 \quad\quad &\{v_1, v_2\}, \{v_3, v_4\}, \{v_5\}, \{v_6\}, \\
&\{e_1, e_4\}, \{e_2, e_3\}, \{e_5, e_6\}, \{e_7, e_8\}, \{e_9\}, \{e_{10}\} \\
\sig_2 \quad\quad &\{v_1\}, \{v_2, v_6\}, \{v_3, v_5\}, \{v_4\}, \\
&\{e_1\}, \{e_2, e_6\}, \{e_3, e_8\}, \{e_4, e_{10}\}, \{e_5, e_7\}, \{e_9\} \\
\sig_3 \quad\quad &\{v_1, v_6\}, \{v_2\}, \{v_3\}, \{v_4, v_5\}, \\
&\{e_1, e_{10}\}, \{e_2, e_7\}, \{e_3, e_5\}, \{e_4\}, \{e_6, e_8\}, \{e_9\} \\
\tau \quad\quad &\{v_1, v_2\}, \{v_2, v_6\}, \{v_3, v_5\}, \\
&\{e_1, e_4, e_{10}\}, \{e_2, e_5, e_8\}, \{e_3, e_7, e_6\}, \{e_9\}
\end{align*}
\caption{Vertex orbits under $S_3$}\label{fig:S3V}
\end{figure}

\begin{figure}[h]
\begin{align*}
\sig_1 \quad\quad &\{\ell_{1,5}, \ell_{4,6}\}, \{\ell_{1,7}, \ell_{4,8}\}, \{\ell_{1,9}, \ell_{4,9}\}, \{\ell_{2,6}, \ell_{3,5}\}, \\
&\{\ell_{2,7}, \ell_{3,8}\}, \{\ell_{2,10}, \ell_{3,10}\}, \{\textcolor{blue}{\ell_{5,6}}\}, \{\textcolor{blue}{\ell_{7,8}}\}, \{\textcolor{red}{\ell_{9,10}}\} \\
\sig_2 \quad\quad &\{\ell_{1,5}, \ell_{1,7}\}, \{\textcolor{red}{\ell_{1,9}}\}, \{\textcolor{blue}{\ell_{2,6}}\}, \{\ell_{2,7}, \ell_{5,6}\}, \{\ell_{2,10}, \ell_{4,6}\}, \\
&\{\ell_{3,5}, \ell_{7,8}\}, \{\textcolor{blue}{\ell_{3,8}}\}, \{\ell_{3,10}, \ell_{4,8}\}, \{\ell_{4,9}, \ell_{9,10}\} \\
\sig_3 \quad\quad &\{\ell_{1,5}, \ell_{3,10}\}, \{\ell_{1,7}, \ell_{2,10}\}, \{\ell_{1,9}, \ell_{9,10}\}, \{\ell_{2,6}, \ell_{7,8}\}, \\
&\{\textcolor{blue}{\ell_{2,7}}\}, \{\textcolor{blue}{\ell_{3,5}}\}, \{\ell_{3,8}, \ell_{5,6}\}, \{\ell_{4,6}, \ell_{4,8}\}, \{\textcolor{red}{\ell_{4,9}}\} \\
\tau \quad\quad &\{\ell_{1,5}, \ell_{4,8}, \ell_{2, 10}\}, \{\ell_{1,7}, \ell_{4,6}, \ell_{3, 10}\}, \{\ell_{1,9}, \ell_{4,9}, \ell_{9, 10}\}, \\
&\{\ell_{2,6}, \ell_{3,5}, \ell_{7, 8}\}, \{\ell_{2,7}, \ell_{5,6}, \ell_{3, 8}\}
\end{align*}
\caption{Edge orbits under $S_3$. \textcolor{blue}{Blue} indicates nontrivial $\bbZ / 2$ action and \textcolor{red}{Red} indicates trivial action.}\label{fig:S3E}
\end{figure}

\begin{figure}[h]
\centering
\begin{tikzpicture}
\draw[thick, domain=-1:6] plot (\x,{0}) node [label={[label distance=0.05cm]0:{$\calL_1$}}] {};
\draw[thick, domain=-0.5:3] plot (\x,{sqrt(3)*\x}) node [label={[label distance=0.05cm]60:{$\calL_1$}}] {};
\draw[thick, domain=2:5.5] plot (\x,{-sqrt(3)*\x+5*sqrt(3)});
\node at (2,5.196) [label={[label distance=0.05cm]120:{$\calL_1$}}] {};
\draw[olive, fill=olive] (1,0) circle (0.1cm) node [label={[label distance=0.05cm]270:{$\calQ_1$}}] {};;
\draw[orange, fill=orange] (2,0) circle (0.1cm) node [label={[label distance=0.05cm]270:{$\calP_1$}}] {};
\draw[red, fill=red] (3,0) circle (0.1cm) node [label={[label distance=0.05cm]270:{$\calQ_2$}}] {};
\draw[violet, fill=violet] (4,0) circle (0.1cm) node [label={[label distance=0.05cm]270:{$\calP_2$}}] {};
\draw[violet, fill=violet] (0.5,0.866) circle (0.1cm) node [label={[label distance=0.05cm]120:{$\calP_2$}}] {};
\draw[olive, fill=olive] (1,1.732) circle (0.1cm) node [label={[label distance=0.05cm]120:{$\calQ_1$}}] {};
\draw[orange, fill=orange] (1.5,2.598) circle (0.1cm) node [label={[label distance=0.05cm]120:{$\calP_1$}}] {};
\draw[red, fill=red] (2,3.464) circle (0.1cm) node [label={[label distance=0.05cm]120:{$\calQ_2$}}] {};
\draw[orange, fill=orange] (4.5,0.866) circle (0.1cm) node [label={[label distance=0.05cm]60:{$\calP_1$}}] {};
\draw[red, fill=red] (4,1.732) circle (0.1cm) node [label={[label distance=0.05cm]60:{$\calQ_2$}}] {};
\draw[violet, fill=violet] (3.5,2.598) circle (0.1cm) node [label={[label distance=0.05cm]60:{$\calP_2$}}] {};
\draw[olive, fill=olive] (3,3.464) circle (0.1cm) node [label={[label distance=0.05cm]60:{$\calQ_1$}}] {};
\end{tikzpicture}
\caption{The cusp triangle orbit of size two}\label{fig:Triangle2}
\end{figure}
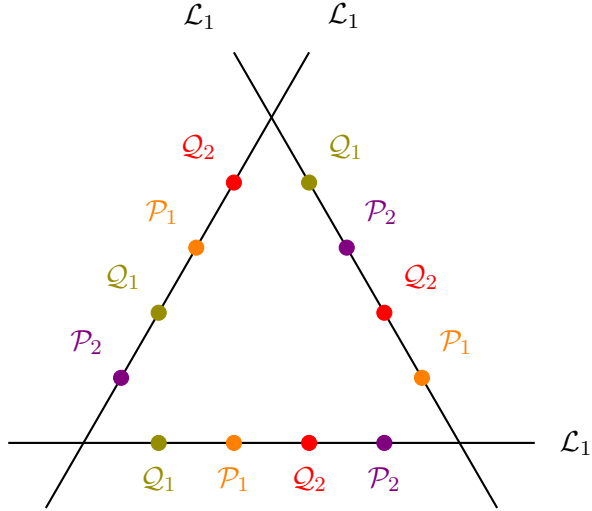

The six icosahedron vertex points $\{v_i\}$ are partitioned into two $S_3$ orbits with three elements:
\begin{align*}
\calP_1 &= \{v_1, v_2, v_6\} \\
\calP_2 &= \{v_3, v_4, v_5\}
\end{align*}
The dodecahedron vertex points $\{e_i\}$ are partitioned as:
\begin{align*}
\calQ_1 &= \{e_1, e_4, e_{10}\} \\
\calQ_2 &= \{e_2, e_3, e_5, e_6, e_7, e_8\} \\
\calQ_3 &= \{e_9\}
\end{align*}
Recall that the lines on $\bbP^2$ associated with the edges of the icosahedron and dodecahedron are the lines arising from the smooth compactification of $X[2]$. These are naturally partitioned into two orbits of triangles, one of size three and one of size two. If $\ell_{i,j}$ is the line containing $e_i$ and $e_j$, we have orbits:
\begin{align*}
\calL_1 &= \{\ell_{1,5}, \ell_{1,7}, \ell_{2,10}, \ell_{3,10}, \ell_{4,6}, \ell_{4,8}\} \\
\calL_2 &= \{\ell_{1,9}, \ell_{4,9}, \ell_{9,10}\} \\
\calL_3 &= \{\ell_{2,6}, \ell_{3,5}, \ell_{7,8}\} \\
\calL_4 &= \{\ell_{2,7}, \ell_{3,8}, \ell_{5,6}\}
\end{align*}
This allows us to visualize each $S_3$ orbit of triangles as in Figure~\ref{fig:Triangle1} and Figure~\ref{fig:Triangle2}.

\begin{figure}[h]
\centering
\begin{tikzpicture}
\draw[thick, domain=-1:6] plot (\x,{0}) node [label={[label distance=0.05cm]0:{$\calL_2$}}] {};
\draw[thick, domain=-0.5:3] plot (\x,{sqrt(3)*\x}) node [label={[label distance=0.05cm]60:{$\calL_3$}}] {};
\draw[thick, domain=2:5.5] plot (\x,{-sqrt(3)*\x+5*sqrt(3)});
\node at (2,5.196) [label={[label distance=0.05cm]120:{$\calL_4$}}] {};
\draw[olive, fill=olive] (1,0) circle (0.1cm) node [label={[label distance=0.05cm]270:{$\calQ_1$}}] {};;
\draw[orange, fill=orange] (2,0) circle (0.1cm) node [label={[label distance=0.05cm]270:{$\calP_1$}}] {};
\draw[violet, fill=violet] (3,0) circle (0.1cm) node [label={[label distance=0.05cm]270:{$\calP_2$}}] {};
\draw[blue, fill=blue] (4,0) circle (0.1cm) node [label={[label distance=0.05cm]270:{$\calQ_3$}}] {};
\draw[red, fill=red] (0.5,0.866) circle (0.1cm) node [label={[label distance=0.05cm]120:{$\calQ_2$}}] {};
\draw[orange, fill=orange] (1,1.732) circle (0.1cm) node [label={[label distance=0.05cm]120:{$\calP_1$}}] {};
\draw[red, fill=red] (1.5,2.598) circle (0.1cm) node [label={[label distance=0.05cm]120:{$\calQ_2$}}] {};
\draw[orange, fill=orange] (2,3.464) circle (0.1cm) node [label={[label distance=0.05cm]120:{$\calP_1$}}] {};
\draw[red, fill=red] (4.5,0.866) circle (0.1cm) node [label={[label distance=0.05cm]60:{$\calQ_2$}}] {};
\draw[violet, fill=violet] (4,1.732) circle (0.1cm) node [label={[label distance=0.05cm]60:{$\calP_2$}}] {};
\draw[red, fill=red] (3.5,2.598) circle (0.1cm) node [label={[label distance=0.05cm]60:{$\calQ_2$}}] {};
\draw[violet, fill=violet] (3,3.464) circle (0.1cm) node [label={[label distance=0.05cm]60:{$\calP_2$}}] {};
\end{tikzpicture}
\caption{The cusp triangle orbit of size three}\label{fig:Triangle1}
\end{figure}
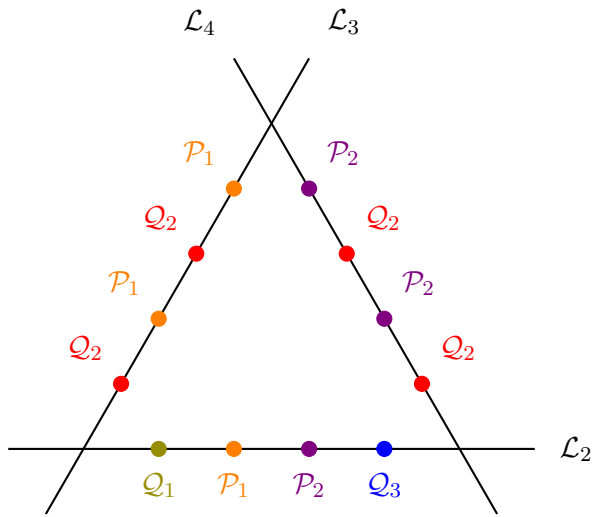

%%%%%%%%%%%%%%%%%%%%
\subsection{The monodromy plane curve}\label{ssec:MPC}
%%%%%%%%%%%%%%%%%%%%

Since $\PSL_2(\bbF_2)$ preserves the image of the lifted period map $\wt{\Pi}^\circ$, we conclude that $\wt{\Pi}^\circ(\wt{\calB}^\circ)$ is an $S_3$-invariant curve of genus one on $X[2]$ with $18$ punctures. Consider the completion
\[
\wt{\Pi} : \wt{\calB} \to Y[2]
\]
and let $\calC = (\sig \circ \wt{\Pi})(\wt{\calB})$ be the image curve on $\bbP^2$. We call $\calC$ the \emph{(Wiman--Edge) monodromy plane curve}. Given the maps and group actions described above, to completely determine the image in $X^\circ$ of the period map for the Wiman--Edge monodromy, it suffices to describe the monodromy plane curve. Proposition~\ref{prop:Bcover} implies the following.

%%%%%%%%%%%%%%%%%%%%
\begin{cor}\label{cor:Bcover}
The monodromy plane curve $\calC$ is a (possibly singular) plane curve of genus one.
\end{cor}
%%%%%%%%%%%%%%%%%%%%

%%%%%%%%%%%%%%%%%%%%
\subsubsection{Intersection of the monodromy plane curve with $\ell_{i,j}$}\label{sssec:MPCcusps}
%%%%%%%%%%%%%%%%%%%%

We now study the intersections of $\calC$ with each line $\ell_{i,j}$. Since $\calC$ is $S_3$-invariant, this only depends on how $\calC$ meets a representative for each orbit. Notice that points on $\calC \cap \ell_{i,j}$ associated with a vertex of the icosahedron or dodecahedron lift to points on $\wt{\Pi}^\circ(\wt{\calB}^\circ)$, not punctures. Conversely, points on $\calC \cap \ell_{i,j}$ that are not a vertex of either the icosahedron or the dodecahedron lift to the cusp divisor for $Y[2]$, hence they are associated with the punctures of $\wt{\Pi}^\circ(\wt{\calB}^\circ)$; we call these \emph{puncture points} on $\ell_{i,j}$. The following should be compared with the description in \cite[\S 5.2]{FL1} of where the punctures of $\calB^\circ$ go in $X^\circ$.

%%%%%%%%%%%%%%%%%%%%
\begin{lem}\label{lem:MPCcusps1}{\ }
\begin{enum}

\item There are either one or two puncture points on each $\ell_{i,j} \in \calL_1$.

\item Each line in $\calL_2$ contains exactly three or four puncture points.

\item Puncture points on lines in $\calL_3$ or $\calL_4$ are fixed points for the $\bbZ / 2$ action on the line, where $\bbZ / 2$ is the stabilizer in $S_3$ of the relevant triangle of lines as in Figure~\ref{fig:Triangle1}. In particular, there are at most two puncture points on any such line.

\end{enum}
The number of puncture points on any line $\ell_{i,j}$ is an invariant of its $S_3$-orbit.
\end{lem}
%%%%%%%%%%%%%%%%%%%%

%%%%%%%%%%%%%%%%%%%%
\begin{pf}
The last statement in the lemma is obvious.

The orbit of $\calL_1$ is associated with the cusp of $(\frakh \times \frakh) / \SL_2(\calO_o)$ denoted $\infty_X$ in \cite[\S 3.3]{FL1}. Indeed, this is the orbit consisting of exactly two triangles of lines on $Y[2]$. In terms of the action of $\Gam$ on $\bbP^1(k)$, $\infty_X$ is associated with the lines $[X : 1]$ and $[1 : X]$ in $\bbP^1(k)$, which represent distinct cusps of $\SL_2(\calO)[2]$ but the same cusp of $\SL_2(\calO_o)$. Since $\gam_\al \gam_\beta \gam_{\al^\prime} \gam_{\beta^\prime}$ stabilizes $[X : 1]$, we see that the six relevant punctures of $\wt{\Pi}^\circ(\wt{\calB}^\circ)$ are given by intersections of $\calC$ with the lines in $\calL_1$.

Similarly, the orbits $\calL_2$, $\calL_3$, and $\calL_4$ are associated with the cusp called $\infty_0$ in \cite[\S 3.3]{FL1}, which is the orbit consisting of exactly three triangles of lines on $Y[2]$. With respect to the action of $\Gam$ on $\bbP^1(k)$, this cusp is associated with the lines $[0 : 1]$, $[1 : 0]$, $[1 : 1]$, and $[-1 : 1]$ in $\bbP^1(k)$, which all represent distinct cusps of $\SL_2(\calO)[2]$ but the same cusp of $\SL_2(\calO_o)$. Since $\gam_\al$ fixes $[1 : 0]$, $\gam_{\al^\prime}$ fixes $[0 : 1]$, and $\gam_\beta$, $\gam_{\beta^\prime}$ both fix $[-1 : 1]$, the twelve relevant punctures of $\wt{\Pi}^\circ(\wt{\calB}^\circ)$ are given by intersections of $\calC$ with the lines in $\calL_2$, $\calL_3$, and $\calL_4$.

Returning to $\calL_1$, the nature of the $S_3$ action implies that each of the two triangles of lines represented by Figure~\ref{fig:Triangle2} determines exactly three punctures of $\wt{\Pi}^\circ(\wt{\calB}^\circ)$. Invariance under the $\bbZ / 3$ action on each individual triangle implies that there are two possibilities: either each line contains a unique puncture point, and that point lies on no other line in the triangle, or the puncture points on the triangle are the three intersection points between the lines, hence each line contains exactly two puncture points. This proves that $\calC$ meets each $\ell_{i,j} \in \calL_1$ as claimed in the statement of the lemma.

Now, consider the triangle associated with lines in $\calL_2, \calL_3, \calL_4$. There are exactly four points in the intersection of $\calC$ with the triangle that are associated with punctures, and these four points are fixed points for the action of $\bbZ / 2$. Since $\bbZ / 2$ acts trivially on the lines in $\calL_2$ and as $-1$ on the curves in $\calL_3$ and $\calL_4$ (under the appropriate identification with $\bbP^1$), puncture points must arise from points on the line in the orbit $\calL_2$ along with possibly the point at the intersection of the lines in $\calL_3$ and $\calL_4$. This leaves only the possibilities in the statement of the lemma.
\end{pf}
%%%%%%%%%%%%%%%%%%%%

%%%%%%%%%%%%%%%%%%%%
\begin{lem}\label{lem:HitQ}
The monodromy plane curve $\calC$ meets one of the $S_3$ orbits $\calQ_j$ associated with the Eckardt points of the Clebsch cubic.
\end{lem}
%%%%%%%%%%%%%%%%%%%%

%%%%%%%%%%%%%%%%%%%%
\begin{pf}
This follows from the description in \cite[\S 5.2]{FL1} of the involution $\iota$, which implies that $\calC$ must meet one of the lifts of $\frakh/\SL_2(\bbZ)$ to $X[2]$. As mentioned above, these lifts map precisely to the points on $\bbP^2$ in the orbits $\calQ_i$ for each $i = 1, 2, 3$.
\end{pf}
%%%%%%%%%%%%%%%%%%%%

%%%%%%%%%%%%%%%%%%%%
\subsubsection{Smooth monodromy plane curve}\label{sssec:MPCsmooth}
%%%%%%%%%%%%%%%%%%%%

In this section, we study the possibility that the monodromy plane curve is smooth. Mathematica code to assist the reader in verifying various assertions that follow is available from the author's webpage \cite{code2}. Since $\calC$ is smooth of genus one, it is a cubic plane curve by the genus-degree formula. We first use $S_3$-invariance to restrict the possibilities for the equation of the curve.

%%%%%%%%%%%%%%%%%%%%
\begin{lem}\label{lem:InvEq}
Suppose that $\calC$ is a smooth plane cubic curve that is invariant under the action of $S_3$ described in \S\ref{sssec:Action} with coordinates $[z_1 : z_2 : z_3]$ on $\bbP^2$. Then $\calC$ is the vanishing set of the homogeneous cubic
\begin{align*}
F\left([z_1 : z_2 : z_3]\right) = a_1 z_1^3 &+ a_2 z_2^3 + a_3 z_2 z_3^2 \\
&+ \Big((2 X - 2) a_1 - X a_2 - (X - 1) a_3\Big) z_1^2 z_2\\ &+ \Big((2 - X) a_1+ a_2 + (X - 1) a_3\Big) z_1 z_2^2 \\ &+\Big((3 X - 5) a_1 + (X + 1) a_2 + (2 - X) a_3\Big) z_1 z_3^2 
\end{align*}
for some $[a_1 : a_2 : a_3] \in \bbP^2$.
\end{lem}
%%%%%%%%%%%%%%%%%%%%

%%%%%%%%%%%%%%%%%%%%
\begin{pf}
Using invariance of $F$ under $\sig_1$, we see that either $F \circ \sig_1 = F$ and every term of $F$ of the form $z_1^{j_1} z_2^{j_2} z_3^{j_3}$ has $j_1 + j_2$ even, or $F \circ \sig_1 = - F$ and every term of $F$ of the form $z_1^{j_1} z_2^{j_2} z_3^{j_3}$ has $j_1 + j_2$ odd. In the first case, one can use a computer algebra program like Mathematica to see that invariance of $F$ under $\sig_2$ implies that $F$ is the singular homogeneous cubic
\[
z_3 \Big((X-2) z_1^2 + 2 z_1 z_2 - (X+1) z_2^2 + z_3^2 \Big).
\]
Under similar analysis in the second case, one sees that $F \circ \sig_2 = - F$ (since $\sig_1$ and $\sig_2$ are conjugate in $S_3$) and furthermore that any such $F$ invariant under $\sig_1$ and $\sig_2$ must then be of the form given in the statement of the lemma. Since $\sig_1$ and $\sig_2$ generate $S_3$, this proves the lemma.
\end{pf}
%%%%%%%%%%%%%%%%%%%%

Note that some of the homogeneous cubics of the form given in the statement of Lemma~\ref{lem:InvEq} are singular. We will not need the finer distinction of which coefficients give a smooth cubic curve.

Bezout's Theorem implies that $\calC$ meets every $\ell_{i,j}$ in two or three points, where the case of two intersection points occurs if and only if $\ell_{i,j}$ is tangent to $\calC$. We begin by completely determining $\calC$ in the case where no $\ell_{i,j}$ is tangent to $\calC$.

%%%%%%%%%%%%%%%%%%%%
\begin{prop}\label{prop:Image1}
Suppose that the monodromy plane curve $\calC$ is smooth. Then $\calC$ is the plane cubic curve with one of the following equations:
\begin{align*}
z_1^3 + (13 - 8 X) z_2^3 &+ (3 - X) z_1^2 z_2 + (18 - 11 X) z_1 z_2^2 \\
&+ (3 X - 5) z_1 z_3^2 + (X - 2) z_2 z_3^2 \\
\\
z_1^3 + (5 - 3 X) z_2^3 &+ X z_1^2 z_2 + (4 - 5 X) z_1 z_2^2 \\ &+ (3 X - 5) z_1 z_3^2 - z_2 z_3^2
\end{align*}
Equivalently, $\calC$ is either
\begin{enum}

\item the unique $S_3$-invariant smooth plane cubic tangent to the lines in $\calL_1$ at the points in $\calQ_2$, or

\item the unique $S_3$-invariant smooth plane cubic containing the orbit $\calQ_2$ (and no other $\calP_i$ or $\calQ_i$) and the intersection points between pairs of lines in the orbit $\calL_1$ represented by the vertices of the triangle in Figure~\ref{fig:Triangle2}.

\end{enum}
\end{prop}
%%%%%%%%%%%%%%%%%%%%

%%%%%%%%%%%%%%%%%%%%
\begin{pf}
Invariance under $S_3$ implies that $\calC$ meets one vertex of the triangle in Figure~\ref{fig:Triangle2} if and only if it meets all the vertices. We break the proof into two cases, accordingly.

\medskip

\noindent
\textbf{Case 1:} No vertex in the triangle from Figure~\ref{fig:Triangle2} is on $\calC$.

\medskip

Lemma~\ref{lem:MPCcusps1} then implies that each line in $\calL_1$ contains exactly one puncture point on $\calC$. We first assume that $\calC$ intersects each line in $\calL_1$ transversally, hence in exactly three points by Bezout's theorem. Then $\calC$ meets exactly two of the orbits $\calP_1$, $\calP_2$, $\calQ_1$, $\calQ_2$. However, $\calC$ cannot simultaneously meet $\calP_1$ and $\calQ_2$ or $\calP_2$ and $\calQ_2$, since otherwise it would meet lines in $\calL_3$ or $\calL_4$ in four points, respectively. This leaves four cases.

\medskip

\noindent
\textbf{Case 1(a):} $\calP_1 \cup \calP_2 \subset \calC$

\medskip

Since $\calC$ does not contain $\calQ_1$ or $\calQ_2$, Lemma~\ref{lem:HitQ} implies that $\calQ_3 \subset \calC$. One checks with a computer algebra program like Mathematica that there is no cubic curve with the form given in Lemma~\ref{lem:InvEq} containing all of $\calP_1 \cup \calP_2 \cup \calQ_3$. Thus this case is impossible.

\medskip

\noindent
\textbf{Case 1(b):} $\calP_1 \cup \calQ_1 \subset \calC$

\medskip

The only cubic as in Lemma~\ref{lem:InvEq} vanishing on $\calP_1 \cup \calQ_1$ is:
\[
z_1^3 - z_2^3 + (X + 1) z_1^2 z_2 + (X - 2) z_1 z_2^2 - (X + 1) z_1 z_3^2 + (2 - X) z_2 z_3^2
\]
However, this cubic also vanishes at $\calP_2$, contradicting our observation that $\calC$ meets exactly two of the orbits $\calP_1$, $\calP_2$, $\calQ_1$, $\calQ_2$. This case is therefore eliminated.

\medskip

\noindent
\textbf{Case 1(c):} $\calP_2 \cup \calQ_1 \subset \calC$

\medskip

The only possibility in this case is the cubic considered in Case 1(b), which eliminates this case as well.

\medskip

\noindent
\textbf{Case 1(d):} $\calQ_1 \cup \calQ_2 \subset \calC$

\medskip

Now, a computer algebra computation shows that the only possibility is the cubic:
\[
(1+X) z_1^3 + z_2^2 + (1 + 3X) z_1^2 z_2 + (1 - 2X) z_1 z_2^2 + (X - 2) z_1 z_3^2 - (2 + 3X) z_2 z_3^2
\]
We rule this case out using the triangle of curves in Figure~\ref{fig:Triangle1}. The associated puncture points on the cubic curve are at the intersection of the lines in $\calL_3$ and $\calL_4$ and two points on the line in $\calL_2$ that do not lie on the other two lines. This implies that the triangle only determines three puncture points, but it must determine four. Indeed, the three triangles making up this $S_3$-orbit of must account for the $12$ punctures associated with $\gam_\al$, $\gam_\beta$, $\gam_{\al^\prime}$, and $\gam_{\beta^\prime}$ (cf.~the proof of Lemma~\ref{lem:MPCcusps1}). This contradiction eliminates this case.

\medskip

We now consider the possibility that $\calC$ is tangent to a line in $\calL_1$. Here $\calC$ meets each line in $\calL_1$ in exactly two points, hence $\calC$ meets exactly one of the orbits $\calP_1$, $\calP_2$, $\calQ_1$, $\calQ_2$.

\medskip

\noindent
\textbf{Case 1(e):} The lines in $\calL_1$ are tangent to $\calC$ at either $\calP_1$, $\calP_2$, or $\calQ_1$.

\medskip

For each of these three possibilities, there is no $S_3$-invariant cubic passing through the relevant orbit that is tangent to $\calL_1$ at the specified points. Therefore, this case is not possible.

\medskip

\noindent
\textbf{Case 1(f):} The lines in $\calL_1$ are tangent to $\calC$ at $\calQ_2$.

\medskip

The only possibility is the first cubic in the statement of the proposition. Then $\calC$ meets the triangle in Figure~\ref{fig:Triangle1} in four puncture points, namely at the intersection of the lines in $\calL_3$ and $\calL_4$ along with three points on the line in $\calL_2$. Therefore this case is a possibility.

\medskip

\noindent
\textbf{Case 1(g):} The lines in $\calL_1$ are tangent to $\calC$ at a puncture point.

\medskip

The cubic must also meet $\calL_1$ in exactly one of $\calP_i$ or $\calQ_i$ ($i=1,2$). For $\calP_1$ and $\calP_2$, the cubic must also meet $\calQ_3$ by Lemma~\ref{lem:HitQ}. For $i=1,2$, the unique cubics through $\calP_i$ and $\calQ_3$ are:
\begin{align*}
&z_1 \big(z_1^2 + (2 X - 2) z_1 z_2 + (2 - X) z_2^2 + (3 X - 5) z_3^2 \big) \\
&z_2 \big(z_2^2 + (X + 1) z_1^2 - 2 X z_1 z_2 - (3 X + 2) z_3^2\big)
\end{align*}
These are singular, which rules these cases out.

For $\calQ_1$, there are cubics that satisfy the required properties with respect to the triangle of lines in $\calL_1$. However, these cubics all would give seven puncture points on the triangle of lines in Figure~\ref{fig:Triangle1}. This contradicts Lemma~\ref{lem:MPCcusps1}, and hence rules out this case. There is no cubic with the required form that meets $\calQ_2$ and is tangent to the lines in $\calL_1$ in a point not in $\calQ_2$. This completes the analysis of this case.

\medskip

The final remaining option for Case 1 is:

\medskip

\noindent
\textbf{Case 1(h):} The lines in $\calL_1$ meet $\calC$ in a triple point.

\medskip

The triple point must be a puncture point of $\calC$, so $\calC$ does not meet $\calP_1$, $\calP_2$, $\calQ_1$, or $\calQ_2$. Lemma~\ref{lem:HitQ} implies that $\calC$ meets $\calQ_3$. One checks that there is no cubic of the appropriate form that meets each line in $\calL_1$ in a triple point, and this completes the analysis of Case 1.

\medskip

\noindent
\textbf{Case 2:} Each vertex in the triangle from Figure~\ref{fig:Triangle2} is on $\calC$.

\medskip

Note that a vertex cannot be an intersection point with multiplicity three. Indeed, then each line in $\calL_1$ would meet $\calC$ in two points, one having multiplicity three, and this contradicts Bezout's theorem. We continue as in the first case by first assuming that no line in $\calL_1$ is tangent to $\calC$. A line in $\calL_1$ then contains two puncture points, hence $\calC$ meets exactly one of the orbits $\calP_1$, $\calP_2$, $\calQ_1$, or $\calQ_2$.

\medskip

\noindent
\textbf{Case 2(a):} $\calP_i \subset \calC$, $i = 1$ or $2$

\medskip

Lemma~\ref{lem:HitQ} implies that $\calC$ also meets $\calQ_3$. One checks that there is no such cubic with the form in Lemma~\ref{lem:InvEq}, hence this case is impossible.

\medskip

\noindent
\textbf{Case 2(b):} $\calQ_1 \subset \calC$

\medskip

This case implies that the cubic is:
\[
z_1^3 + (1 + X) z_2^3 + X z_1^2 z_2 - X z_1 z_2^2 + (3 X - 1) z_1 z_3^2 - (1 + 4 X) z_2 z_3^2
\]
This cubic curve meets the triangle in Figure~\ref{fig:Triangle2} in what would be seven puncture points. In particular, it meets the triangle at the intersection of the lines in $\calL_3$ and $\calL_4$, two additional points on each of those lines, and two points on the line in $\calL_2$ (neither of which is $e_9 \in \calQ_3$). This contradicts the fact that this triangle must determine exactly four puncture points, as we saw in Case 1(d). This case is therefore eliminated.

\medskip

\noindent
\textbf{Case 2(c):} $\calQ_2 \subset \calC$

\medskip

The only possibility here is the second cubic given in the statement of the proposition. As in Case 1(h), $\calC$ then meets the triangle in Figure~\ref{fig:Triangle1} in four puncture points, namely at the intersection of the lines in $\calL_3$ and $\calL_4$ along with three points on the line in $\calL_2$. Therefore this case is a possibility.

\medskip

\noindent
\textbf{Case 2(d):} The lines in $\calL_1$ are tangent to $\calC$.

\medskip

In this case, Bezout's theorem implies that $\calC$ meets the triangle of lines in $\calL_1$ in only the three vertex points. Specifically, it suffices to study plane cubics passing through $\ell_{1,5} \cap \ell_{4,8}$ that are tangent to either $\ell_{1,5}$ or $\ell_{4,8}$. These are:
\begin{align*}
&z_1^3 + (X-1) z_2^3 + X z_1^2 z_2 + (X-2) z_1 z_2^2 + (X-1) z_1 z_3^2 + (1- 2 X) z_2 z_3^2 \\
&z_1^3 + (1-X) z_2^3 + X z_1^2 z_2 + (2-X) z_1 z_2^2 + (X-3) z_1 z_3^2 + z_2 z_3^2
\end{align*}
Neither cubic can vanish at any other point on the triangle, hence it cannot vanish on the orbits $\calQ_1$ or $\calQ_2$. Lemma~\ref{lem:HitQ} then implies that it must vanish on $\calQ_3 = \{e_9\}$. Neither cubic is zero at this point, hence we have eliminated the possibility that $\calC$ is tangent to a line in $\calL_1$. This completes the analysis of Case 2, and hence the proof of the proposition.
\end{pf}
%%%%%%%%%%%%%%%%%%%%

%%%%%%%%%%%%%%%%%%%%
\bibliography{Wiman}
%%%%%%%%%%%%%%%%%%%%

%%%%%%%%%%%%%%%%%%%%
\end{document}